\pdfoutput=1
\documentclass[11pt]{amsart}
\usepackage[margin=1in]{geometry}
\usepackage{amsmath,amssymb,amsthm}
\usepackage{enumitem}
\usepackage{xcolor}
\usepackage[hidelinks]{hyperref}
\usepackage{refcount}
\usepackage{tikz}
\usetikzlibrary{arrows.meta,positioning,shapes.geometric}
\DeclareRobustCommand{\symdone}{\tikz[baseline=-.6ex]\node[circle, fill=black, inner sep=0pt, minimum size=1.8mm]{};}
\DeclareRobustCommand{\symdh}{\tikz[baseline=-.6ex]\node[rectangle, fill=black, inner sep=0pt, minimum size=1.6mm]{};}
\DeclareRobustCommand{\symdv}{\tikz[baseline=-.6ex]\node[diamond, fill=black, inner sep=0pt, minimum size=2.2mm]{};}
\DeclareRobustCommand{\symdd}{\tikz[baseline=-.6ex]\node[circle, fill=white, draw=black, line width=.5pt, inner sep=0pt, minimum size=1.8mm]{};}

\setlist{itemsep=2pt,topsep=4pt}

\newtheorem{theorem}{Theorem}[section]
\newtheorem{proposition}[theorem]{Proposition}
\newtheorem{lemma}[theorem]{Lemma}
\newtheorem{corollary}[theorem]{Corollary}
\newtheorem{maintheorem}{Theorem}
\newtheorem{maincorollary}{Corollary}

\theoremstyle{definition}

\newtheorem{example}[theorem]{Example}
\theoremstyle{remark}
\newtheorem{note}[theorem]{Note}

\title{The Galois actions on braided fusion categories}
\author{Andrew Schopieray}
\address{Mathematics Department, SUNY Plattsburgh, Plattsburgh, NY 12901, USA}
\email{ascho013@plattsburgh.edu}

\begin{document}

\begin{abstract}
This note harmonizes two classical topics: minimal modular extensions of braided fusion categories and the Galois action on modular fusion categories. Their common generalizations are Galois-modular extensions, which always exist when a braided fusion category $\mathcal{B}$ is pseudounitary. Each extension induces a Galois action on the simple objects of $\mathcal{B}$. Up to integral equivalence, the Galois-modular extensions of $\mathcal{B}$ form a torsor over an explicit group of homomorphisms from $\operatorname{Gal}(\overline{\mathbb{Q}}/\mathbb{Q})$ to the invertible objects of order dividing two in the symmetric center $\mathcal{E}$ of $\mathcal{B}$. If $\mathcal{E}$ is Tannakian and $\mathcal{B}$ has a minimal modular extension, all minimal modular extensions induce one canonical Galois action. We prove that many facts about the modular data of modular fusion categories extend to pseudounitary premodular fusion categories as well via the Galois action.
\end{abstract}

\maketitle

\null
\vspace{-10mm}
% ======================================================================
\section{A non-technical introduction}
% ======================================================================

\par M.\ M\"uger conjectured \cite[Conjecture 5.2]{MR1990929} that unitary premodular fusion categories could always be extended to modular fusion categories of the smallest possible dimension. Though the conjecture is false \cite{GVR,MR4771215}, it instigated a line of research which continues today. Much of this research concerns determining for which categories the conjecture is true, studying minimal modular extensions when they exist, or discussing the ramifications of existence to topics in mathematical physics.

\par Our approach is quite different. J.\ Plavnik, Z.\ Yu, Q.\ Zhang, and the author proved a result which implies that a premodular category is stable under the classical Galois action of any modular extension with integral centralizer \cite[Theorem 4.1.6]{MR4834521}. Motivated by this, the author proposed a replacement for the minimal modular extension conjecture: every premodular fusion category has such an extension, called a Galois-modular extension \cite{SchopierayBIRS2026}. The many results which follow from the Galois action on modular fusion categories then pass to the far more general class of premodular categories. T.\ Johnson-Freyd has since proved the existence of Galois-modular extensions for pseudounitary categories using the tools of higher category theory \cite{ART_2026__3_3_579_0}.

\par Now equipped with existence, we show that the basic extension theory of T.\ Lan, L.\ Kong, and X.-G.\ Wen \cite{MR3613518} generalizes to modular extensions which need not be minimal. Roughly, any nondegenerate extension of a braided fusion category $\mathcal{B}$ with a fixed centralizer $\mathcal{C}$ can be obtained from another by a known construction involving minimal nondegenerate extensions of the symmetric center $\mathcal{E}$ of $\mathcal{B}$.  When these nondegenerate extensions are Galois-modular extensions, the relationship between the induced Galois actions on $\mathcal{B}$ from such extensions can be described explicitly in terms of one another. This provides a family of Galois actions on any pseudounitary braided fusion category.  Moreover, these induced Galois actions act on the premodular data, and classical results describing this action and the number fields involved follow in turn. Most of these generalizations are identical to the classical case when $\mathcal{E}$ is Tannakian, e.g.\ for modular fusion categories themselves, and involve only slight modifications for non-Tannakian symmetric centers.

\par Following this non-technical introduction, we provide a technical introduction for experts in Section \ref{sec:tech} which includes statements of the main results. Section \ref{sec:prelim} contains the preliminary notions needed for the proofs of these results.  Section \ref{sec:structure} proves the main structure theorem for nondegenerate extensions and Section \ref{sec:galois} restricts to Galois-modular extensions and their induced Galois actions.  Lastly, Section \ref{sec:data} describes how these induced actions interact with premodular data.

% ======================================================================
\section{Technical introduction}\label{sec:tech}
% ======================================================================

\par Galois automorphisms act both on the characters and the set of simple objects $\mathcal{O}(\mathcal{M})$ of modular fusion categories $\mathcal{M}$ \cite{MR1266785}. In particular, the characters of the fusion algebra $K(\mathcal{M})\otimes_\mathbb{Z}\mathbb{C}$ are in bijection with $Y\in\mathcal{O}(\mathcal{M})$ via $\varphi_Y(X)=s_{X,Y}/d_Y$ where $s$ is the unnormalized $S$-matrix of $\mathcal{M}$ and $d_Y$ is the categorical dimension of $Y$. Therefore each $\sigma\in\operatorname{Gal}(\overline{\mathbb{Q}}/\mathbb{Q})$ determines a permutation $\hat\sigma$ of $\mathcal{O}(\mathcal{M})$ by $\sigma\circ\varphi_Y=\varphi_{\hat\sigma(Y)}$. For a general degenerate premodular category $\mathcal{B}$, to our knowledge no such action on $\mathcal{O}(\mathcal{B})$ is described in the literature, since the characters $\varphi_Y$ no longer separate the simple objects. A \emph{Galois-modular extension} of $\mathcal{B}$ is a braided fully faithful functor $\iota\colon\mathcal{B}\to\mathcal{M}$ into a modular fusion category whose centralizer $C_\mathcal{M}(\iota\mathcal{B})$ is integral. By \cite[Theorem 4.1.6]{MR4834521}, $\iota(\mathcal{O}(\mathcal{B}))$ is stable under $\hat\sigma$, so $\iota$ induces an action $\hat\sigma_\iota$ of $\operatorname{Gal}(\overline{\mathbb{Q}}/\mathbb{Q})$ on $\mathcal{O}(\mathcal{B})$. T.\ Johnson-Freyd proved that every pseudounitary braided fusion category has a Galois-modular extension \cite{ART_2026__3_3_579_0}. The induced action depends on $\iota$ and this paper describes precisely how.

\par Most of the data of this Galois action comes from a generalization of results of T.\ Lan, L.\ Kong and X.-G.\ Wen, who proved that the minimal modular extensions of $\mathcal{B}$, those with centralizer equal to the symmetric center $\mathcal{E}$ of $\mathcal{B}$, form a torsor over the group $\mathrm{Mext}(\mathcal{E})$ of minimal modular extensions of $\mathcal{E}$ \cite{MR3613518}. Minimal modular extensions are Galois-modular, because $\mathcal{E}$ is integral, but they need not exist: for $G$-equivariantizations $\mathcal{B}=\mathcal{C}^G$ of a nondegenerate $\mathcal{C}$ with a braided action of a finite group $G$ they amount to $G$-crossed extensions of $\mathcal{C}$, whose existence is obstructed by a class in $H^4(G,\mathbb{C}^\times)$ \cite{MR2677836,GVR}, \cite[Corollary 1.4]{JOY}; they always exist when $\mathcal{E}\simeq\mathrm{sVec}$ \cite{JFR}. For elementary non-examples, large families of braided Tambara--Yamagami categories have no minimal modular extension \cite{MR4771215}, including representation categories of extra-special $2$-groups with carefully selected braidings.

\par The primary construction from Lan--Kong--Wen is left unchanged although we give it a name so that we can refer to it more casually. For nondegenerate extensions $\iota:\mathcal{B}\to\mathcal{M}$ and $j\colon\mathcal{E}_1\to\mathcal{N}$ of a fusion subcategory $\mathcal{E}_1\subseteq\mathcal{E}$, the \emph{Baer sum} $\iota\oplus j$ embeds $\mathcal{B}$ into the category $\mathcal{M}\boxtimes_{\mathcal{E}_1}\mathcal{N}$ of local modules over the connected \'etale algebra $\bigoplus_{V\in\mathcal{O}(\mathcal{E}_1)}V\boxtimes j(V)^\ast$ (refer to Section \ref{sec:structure}). For minimal nondegenerate extensions of $\mathcal{E}\simeq\mathrm{Rep}(G)$ this is literally a Baer sum of extensions of $G$ by categorical $2$-groups (Note \ref{rem:baer}). Write $\mathcal{B}\boxtimes_{\mathcal{E}}\mathcal{C}$ for the relative Deligne product (Section \ref{subsec:relative}). The structure theorem (Theorem~\ref{thm:structure}) for arbitrary nondegenerate extensions then follows, with no minimality assumption: for any braided fusion category $\mathcal{B}$, a nondegenerate extension of $\mathcal{B}$ with centralizer $\mathcal{C}$ is the same as a minimal nondegenerate extension of $\mathcal{B}\boxtimes_\mathcal{E}\mathcal{C}$, and any two pseudounitary nondegenerate extensions of $\mathcal{B}$ become equivalent after Baer sums with nondegenerate extensions of $\mathcal{E}$ (Lan--Kong--Wen state their results for unitary categories, but their arguments apply to pseudounitary ones, which is all that we need). The following large swath of vocabulary and definitions is specific to our main result.

\par Assume that $\mathcal{B}$ is pseudounitary, in the sense that the Frobenius--Perron and global dimension coincide. A nondegenerate extension of $\mathcal{B}$ with integral centralizer $\mathcal{C}$ is then automatically pseudounitary by \cite[Proposition 8.24]{MR2183279} and \cite[Proposition 2.4]{ART_2026__3_3_579_0}, hence Galois-modular. So by Theorem~\ref{thm:structure}(a), the Galois-modular extensions of $\mathcal{B}$ are exactly the minimal nondegenerate extensions of the categories $\mathcal{B}\boxtimes_\mathcal{E}\mathcal{C}$ with $\mathcal{C}$ integral and $\mathcal{C}'=\mathcal{E}$. For the induced Galois action, define $u_\iota(\sigma):=\hat\sigma_\iota(\mathbf{1})$, write $\mathcal{E}_{\mathrm{pt}}[2]$ for the group of invertible objects of $\mathcal{E}$ of order dividing $2$, and let $\kappa$ be the quadratic character of $\mathbb{Q}(\sqrt2)$, so that $\kappa(\sigma_\ell)=\bigl(\tfrac{2}{\ell}\bigr)$ with $\sigma_\ell$ as in Section \ref{subsec:galois-action}; for an integer $\ell$, $\sigma_\ell$ means any $\sigma$ acting by $\zeta\mapsto\zeta^\ell$ on the cyclotomic field at hand. Let $\mathcal{A}(\mathcal{B})$ be the set of actions on $\mathcal{O}(\mathcal{B})$ induced by Galois-modular extensions of $\mathcal{B}$. A real multiquadratic field is a compositum of finitely many real quadratic fields, $\mathbb{Q}$ included. Call a nondegenerate extension $z\colon\mathcal{E}\to\mathcal{P}$ \emph{integral} if the category $\mathcal{P}$ itself is integral, hence it is Galois-modular. By Lemma~\ref{lem:symmetric}(c), the integral extensions of $\mathcal{E}$ are exactly the Galois-modular ones with $u_z\equiv\mathbf{1}$, so by part (a) of Theorem~\ref{thm:main-intro} below, Baer sums with them do not change the induced Galois action. For this reason, call Galois-modular extensions $\iota,\iota'$ of $\mathcal{B}$ \emph{integrally equivalent} if $\iota\oplus z\simeq\iota'\oplus z'$ for some integral extensions $z,z'$ of $\mathcal{E}$.

\vspace{5mm}

\begin{maintheorem}\label{thm:main-intro}
Let $\mathcal{B}$ be a pseudounitary braided fusion category with symmetric center $\mathcal{E}$.
\begin{enumerate}[label=(\alph*)]
\item If $\iota,j$ are Galois-modular extensions of $\mathcal{B}$ and a fusion subcategory $\mathcal{E}_1\subseteq\mathcal{E}$, respectively, then $\iota\oplus j$ is Galois-modular, and $\hat\sigma_{\iota\oplus j}(X)=\hat\sigma_\iota(X)\otimes u_j(\sigma)$ for all $X\in\mathcal{O}(\mathcal{B})$.
\item Under Baer sums along $\mathcal{E}$, the integral equivalence classes of Galois-modular extensions of $\mathcal{E}$ form an abelian group. The map $j\mapsto u_j$ is an isomorphism from it onto the set $\mathcal{U}(\mathcal{E})$ of all $u_j$, a group under pointwise product. The integral equivalence classes of Galois-modular extensions of $\mathcal{B}$ form a torsor over this group, and $\iota\mapsto\hat\sigma_\iota$ is an equivariant bijection onto $\mathcal{A}(\mathcal{B})$. In particular, $\mathcal{A}(\mathcal{B})$ is a $\mathcal{U}(\mathcal{E})$-torsor under $(\hat\sigma_\iota,u)\mapsto\hat\sigma_\iota\otimes u$, and two Galois-modular extensions of $\mathcal{B}$ induce the same action if and only if they are integrally equivalent.
\item $\mathcal{U}(\mathcal{E})$ is the group of homomorphisms $u\colon\operatorname{Gal}(\overline{\mathbb{Q}}/\mathbb{Q})\to\mathcal{E}_{\mathrm{pt}}[2]$ that factor through $\operatorname{Gal}(K/\mathbb{Q})$ for some real multiquadratic field $K$ and satisfy $\theta\circ u\in\{1,\kappa\}$, where $\theta$ is the twist.
\end{enumerate}
\end{maintheorem}

The payoff is that the induced Galois action of a Galois-modular extension is determined by $\mathcal{B}$ up to a twist by an element of $\mathcal{U}(\mathcal{E})$ (refer to Figure \ref{fig:torsor} for an example); for modular $\mathcal{B}$ it is unique. Conversely, the action determines the extension up to Baer sums with integral extensions of $\mathcal{E}$, which do not change it. Two consequences are immediate.

\begin{maincorollary}
Let $G$ be a finite group. Every induced Galois action on $\mathrm{Rep}(G)$ has the form $X\mapsto u(\sigma)\otimes(\sigma\circ X)$. Here $\sigma\circ X$ is the Galois conjugate representation, and $u\colon\operatorname{Gal}(\overline{\mathbb{Q}}/\mathbb{Q})\to\operatorname{Hom}(G,\{\pm1\})$ is a homomorphism that factors through $\operatorname{Gal}(K/\mathbb{Q})$ for some real multiquadratic field $K$. Every such $u$ occurs.
\end{maincorollary}

\begin{maincorollary}
Let $\mathcal{B}$ be a pseudounitary braided fusion category with symmetric center $\mathcal{E}$. If $\iota,\iota'$ are Galois-modular extensions of $\mathcal{B}$ with (additive) central charges $c,c'$, then $2(c-c')\in\mathbb{Z}$, and $c-c'\in\mathbb{Z}$ if $\mathcal{E}$ is Tannakian. If $\mathcal{E}$ is Tannakian and $\mathcal{B}$ has minimal modular extensions, then all minimal modular extensions induce the same Galois action on $\mathcal{B}$.
\end{maincorollary}

As an application of this novel Galois action, we prove results analogous to (a portion of) \cite{MR1266785} and \cite[Theorem II]{MR3435813}: we describe how the induced Galois actions act on the normalized premodular data of a premodular fusion category $\mathcal{B}$ (Lemma \ref{cond:lem-unit}) and constrain the cyclotomic fields containing it (Section \ref{sec:conductor}).  Proposition~\ref{cond:prop-reverse} proves, in a sharper form, the case of a conjecture of P.\ Bruillard \cite[Conjecture 3.5]{MR3548123}, on the cyclotomic field containing the categorical dimensions of simple objects of a premodular category, for Tannakian symmetric center; the proof of $N_\Phi\mid N_\theta$ uses only that $\mathcal{E}\simeq\mathrm{Rep}(G)$ has trivial twists, not pseudounitarity or the Galois action. As Bruillard observes, its factor $2$ is necessary in general (Note~\ref{cond:note}(ii)); to our knowledge the non-Tannakian case remains open. For super-modular categories, where $\mathcal{E}\simeq\mathrm{sVec}$, the $S$-matrix factors through the fermionic quotient \cite[Theorem 3.5]{MR3641612}, and a Galois action on its labels, the pairs $\{X,f\otimes X\}$, together with the number fields generated by its $S$-matrix, was studied in \cite[Section 2.5]{BPRZ2021} and used in the classification results of \cite{BGNPRW2020,BPRZ2021}. Lemmas~\ref{cond:lem-orbits} and~\ref{cond:lem-unit} extend this action to every pseudounitary premodular category, with $\mathcal{E}$-orbits in place of these pairs, and Theorem~\ref{thm:main-intro} lifts it to the simple objects themselves: by part (c), $\mathcal{U}(\mathrm{sVec})$ has order two, so a pseudounitary super-modular category has exactly two induced Galois actions, both minimal \cite{JFR}, and they differ by tensoring with $f$ exactly when $\sigma(\sqrt2)=-\sqrt2$.

\subsection*{Global conventions}
All of our categories are over $\mathbb{C}$. We write $C_\mathcal{M}(\mathcal{D})$ for the relative centralizer of $\mathcal{D}$ in $\mathcal{M}$ and $\mathcal{D}'=C_\mathcal{D}(\mathcal{D})$ for the symmetric center; $\mathcal{M}$ is called nondegenerate if $\mathcal{M}'=\mathrm{Vec}$. A nondegenerate extension of $\mathcal{B}$ is a braided fully faithful functor $\iota\colon\mathcal{B}\to\mathcal{M}$ with $\mathcal{M}$ nondegenerate. Two extensions $\iota,\iota'$ are equivalent if there is a braided equivalence $F$ from the target of $\iota$ to that of $\iota'$ with $F\iota\cong\iota'$ as braided tensor functors. $\mathrm{Mext}(\mathcal{D})$ is the set of equivalence classes of nondegenerate extensions $\mathcal{M}\supseteq\mathcal{D}$ with $C_\mathcal{M}(\mathcal{D})=\mathcal{D}'$. We use repeatedly that $\mathrm{FPdim}(\mathcal{D})\mathrm{FPdim}(C_\mathcal{M}(\mathcal{D}))=\mathrm{FPdim}(\mathcal{M})$ and $C_\mathcal{M}(C_\mathcal{M}(\mathcal{D}))=\mathcal{D}$ for a fusion subcategory $\mathcal{D}$ of a nondegenerate $\mathcal{M}$ \cite[Theorems 3.10 and 3.14]{MR2609644}. For invertible $g$ and simple $Y$, $b(g,Y)$ denotes the scalar by which $c_{Y,g}c_{g,Y}$ acts; $b(g,Y)=1$ if and only if $g$ and $Y$ centralize each other.

Sections~\ref{subsec:local}, \ref{subsec:relative} and~\ref{sec:structure} need no further assumptions, except where pseudounitarity is stated explicitly. In Section~\ref{subsec:galois-action} and from Section~\ref{sec:galois} on, braided fusion categories are pseudounitary with their canonical spherical structures, so $d_X=\mathrm{FPdim}(X)>0$, and a modular fusion category is a pseudounitary nondegenerate braided one. A nondegenerate extension of a pseudounitary category with pseudounitary centralizer is pseudounitary \cite[Proposition 2.4]{ART_2026__3_3_579_0}; we then call it a \emph{modular extension}. This applies to the elements of $\mathrm{Mext}(\mathcal{D})$, and to extensions with integral centralizer \cite[Proposition 8.24]{MR2183279}. We write $s_{X,Y}$ for the unnormalized $S$-matrix and $\theta$ for the twist; then $\varphi_Y(g)=b(g,Y)$.

% ======================================================================
\section{Preliminaries}\label{sec:prelim}
% ======================================================================

\subsection{Local modules}\label{subsec:local}

Following \cite{MR3039775}, let $\mathcal{P}$ be a nondegenerate braided fusion category and $A\in\mathcal{P}$ a connected \'etale algebra. Write $\mathcal{P}_A$ for the fusion category of right $A$-modules, $\mathcal{P}_A^0$ for its subcategory of local modules, $C_\mathcal{P}(A)$ for the centralizer of the subcategory generated by the simple summands of $A$, and $F(Y)=Y\otimes A$.

\begin{lemma}\label{lem:trace}
Let $\mathcal{P}$ be pseudounitary, and let $L\in\mathcal{P}_A^0$ be simple, with underlying object $\bigoplus_iY_i$ in $\mathcal{P}$, and let $W\in C_\mathcal{P}(A)$ be simple. Then $\mathcal{P}_A^0$ is pseudounitary, $d_{F(W)}=d_W$ and
\[
s^{\mathcal{P}_A^0}_{L,F(W)}=\frac1{\dim A}\sum_is^{\mathcal{P}}_{Y_i,W}.
\]
Moreover, if $a$ is an invertible summand of $A$, then $b(a,Y_i)=1$ for all $i$.
\end{lemma}

\begin{proof}
Since $\mathcal{P}$ is pseudounitary, $\theta_A=\mathrm{id}_A$ \cite[Lemma 2.2.4]{MR3808050}, and $A$ is simple as an $A$-module because it is connected, so $A$ is a rigid algebra in the sense of \cite{MR1936496} by \cite[Lemma 1.20]{MR1936496}. Hence the dimension of an $A$-module is its dimension in $\mathcal{P}$ divided by $\dim A$ \cite[Theorem 1.18]{MR1936496}. In particular $d_{F(W)}=d_W$, and dimensions in $\mathcal{P}_A^0$ are positive, so $\mathcal{P}_A^0$ is pseudounitary and the spherical structure of \cite{MR1936496} on $\mathcal{P}_A^0$ is the canonical one. As $W$ centralizes $A$, the module $F(W)$ is local, and \cite[Theorem 4.1]{MR1936496} gives $\dim A\cdot s^{\mathcal{P}_A^0}_{L,F(W)}=\sum_{Y\in\mathcal{O}(\mathcal{P})}\dim\operatorname{Hom}_A(F(Y),L)\,s^{\mathcal{P}}_{Y,W}$. By adjunction \cite[Theorem 1.6]{MR1936496}, $\dim\operatorname{Hom}_A(F(Y),L)$ is the multiplicity of $Y$ in $\bigoplus_iY_i$, which gives the formula for $s$. For the last claim, the action $\mu\colon L\otimes A\to L$ restricts to an isomorphism $L\otimes a\to L$. Indeed, since $A$ is connected \'etale, the pairing $A\otimes A\to A\to\mathbf{1}$ is nondegenerate, so the product $a\otimes a^{-1}\to A$ has a nonzero component to $\mathbf{1}$, and by associativity $\mu\circ(\mu|_{L\otimes a}\otimes\mathrm{id}_{a^{-1}})$ is a nonzero multiple of the identity. Locality says $\mu\circ c_{A,L}c_{L,A}=\mu$. Restricted to $L\otimes a$, this gives $c_{a,L}c_{L,a}=\mathrm{id}$.
\end{proof}

%%%%

\subsection{The relative Deligne product}\label{subsec:relative}

Let $\mathcal{E}$ be a symmetric braided fusion category; we have $\mathcal{E}^{\mathrm{rev}}=\mathcal{E}$ in this case. For braided fusion categories $\mathcal{B},\mathcal{C}$ with $\mathcal{E}\subseteq\mathcal{B}'$ and $\mathcal{E}\subseteq\mathcal{C}'$, let
\[
L_\mathcal{E}:=\bigoplus_{V\in\mathcal{O}(\mathcal{E})}V\boxtimes V^\ast\in\mathcal{E}\boxtimes\mathcal{E}\subseteq\mathcal{B}\boxtimes\mathcal{C}.
\]
This is the canonical algebra of the tensor product functor $\mathcal{E}\boxtimes\mathcal{E}\to\mathcal{E}$. It is a connected \'etale algebra. Define $\mathcal{B}\boxtimes_\mathcal{E}\mathcal{C}:=(\mathcal{B}\boxtimes\mathcal{C})_{L_\mathcal{E}}$. Every object of $\mathcal{B}\boxtimes\mathcal{C}$ centralizes $\mathcal{E}\boxtimes\mathcal{E}$, so all modules are local and $\mathcal{B}\boxtimes_\mathcal{E}\mathcal{C}$ is braided. It has Frobenius--Perron dimension $\mathrm{FPdim}(\mathcal{B})\mathrm{FPdim}(\mathcal{C})/\mathrm{FPdim}(\mathcal{E})$, and $X\mapsto F(X\boxtimes\mathbf{1})$ and $Y\mapsto F(\mathbf{1}\boxtimes Y)$ embed $\mathcal{B}$ and $\mathcal{C}$ into it. The operation $\boxtimes_\mathcal{E}$ is associative and commutative, with unit $\mathcal{E}$ \cite{MR3022755,MR3613518}. For nondegenerate extensions $\mathcal{M},\mathcal{N}\supseteq\mathcal{E}$, denote $\mathcal{M}\boxtimes_\mathcal{E}\mathcal{N}:=(\mathcal{M}\boxtimes\mathcal{N})^0_{L_\mathcal{E}}$.

\begin{lemma}\label{lem:join}
Let $\mathcal{B},\mathcal{C}$ be fusion subcategories of a braided fusion category that centralize each other, with $\mathcal{B}\cap\mathcal{C}=\mathcal{E}$ symmetric. Then $X\boxtimes Y\mapsto X\otimes Y$ induces a braided equivalence $\mathcal{B}\boxtimes_{\mathcal{E}}\mathcal{C}\simeq\mathcal{B}\vee\mathcal{C}$.
\end{lemma}

\begin{proof}
Since $\mathcal{B}$ and $\mathcal{C}$ centralize each other, $T(X\boxtimes Y)=X\otimes Y$ is a braided tensor functor, and it induces a functor $\bar T\colon\mathcal{B}\boxtimes_{\mathcal{E}}\mathcal{C}\to\mathcal{B}\vee\mathcal{C}$ through the algebra map $T(L_\mathcal{E})\to\mathbf{1}$. Every simple object of $\mathcal{B}\vee\mathcal{C}$ is a summand of some $X\otimes Y$. For $X,X'\in\mathcal{B}$ and $Y,Y'\in\mathcal{C}$,
\[
\operatorname{Hom}(X\otimes Y,X'\otimes Y')\cong\operatorname{Hom}(X'^*\otimes X,\,Y'\otimes Y^*),
\]
and the common simple constituents of an object of $\mathcal{B}$ and an object of $\mathcal{C}$ lie in $\mathcal{E}$. So this space has dimension $\sum_{V\in\mathcal{O}(\mathcal{E})}\dim\operatorname{Hom}(X,X'\otimes V)\dim\operatorname{Hom}(Y,Y'\otimes V^\ast)$. By adjunction \cite[Theorem 1.6]{MR1936496}, that is $\dim\operatorname{Hom}(F(X\boxtimes Y),F(X'\boxtimes Y'))$. A tensor functor is faithful, so $\bar T$ is fully faithful on free modules, hence an equivalence onto $\mathcal{B}\vee\mathcal{C}$.
\end{proof}

Lastly, we recall three of the major results of Lan--Kong--Wen, \cite[Theorems 4.20, 5.4 and 4.22]{MR3613518}, which are parts (a), (b) and (c) below, respectively. They are stated for unitary categories; the arguments apply to pseudounitary categories \cite[Remark 3.1]{ART_2026__3_3_579_0}; for $\mathrm{Mext}(\mathcal{E})$ see also \cite{Nikshych2022}.

\begin{theorem}\label{thm:LKW}
Let $\mathcal{E}$ be a symmetric braided fusion category.
\begin{enumerate}[label=(\alph*)]
\item $\mathrm{Mext}(\mathcal{E})$ is a group under $\boxtimes_{\mathcal{E}}$, with unit $\mathcal{Z}(\mathcal{E})$.
\item If $\mathcal{D}$ is pseudounitary, $\mathcal{D}'=\mathcal{E}$ and $\mathrm{Mext}(\mathcal{D})\ne\emptyset$, then it is a $\mathrm{Mext}(\mathcal{E})$-torsor under $\boxtimes_{\mathcal{E}}$.
\item If $\mathcal{E}\simeq\mathrm{Rep}(G)$ is Tannakian, $\mathrm{Mext}(\mathcal{E})$ consists of the classes of $\mathcal{E}\subseteq\mathcal{Z}(\mathrm{Vec}_G^\omega)$, and $\omega\mapsto\mathcal{Z}(\mathrm{Vec}_G^\omega)$ is an isomorphism $H^3(G,\mathbb{C}^\times)\cong\mathrm{Mext}(\mathcal{E})$.
\end{enumerate}
\end{theorem}

%%%%%%

\subsection{The Galois action}\label{subsec:galois-action}

In this preliminary subsection, and from Section \ref{sec:galois} on, braided fusion categories are pseudounitary. Every $\sigma\in\operatorname{Gal}(\overline{\mathbb{Q}}/\mathbb{Q})$ acts on roots of unity by $\zeta\mapsto\zeta^\ell$ for a unique $\ell\in\widehat{\mathbb{Z}}^\times$, and we write $\sigma=\sigma_\ell$ when only this action matters. For a modular fusion category $\mathcal{M}$, $\hat\sigma$ is defined by $\sigma\circ\varphi_Y=\varphi_{\hat\sigma(Y)}$ \cite{MR1266785}. Since $s_{X,Y^\ast}=\overline{s_{X,Y}}$ and the values of $\varphi_Y$ lie in a cyclotomic field where any $\sigma$ commutes with complex conjugation, the Galois action commutes with duality. Let $p^\pm=\sum_{X\in\mathcal{O}(\mathcal{M})}\theta_X^{\pm1}d_X^2$. The first multiplicative central charge is $\xi_\mathcal{M}=p^+/\sqrt{\dim(\mathcal{M})}=e^{2\pi ic/8}$, where $c\in\mathbb{Q}/8\mathbb{Z}$ is the (additive) central charge. Let $\gamma$ be a cube root of $\xi_\mathcal{M}$; it is a valid choice of the sixth root $\zeta$ of $p^+/p^-=\xi_\mathcal{M}^2$ in \cite{MR3435813}. Since $\ell^2\equiv1\pmod3$, the quotient $\gamma/\sigma^2(\gamma)$ does not depend on the choice of $\gamma$. If $g$ is invertible, then symmetry of $s$ and multiplicativity of $\varphi_Y$ give $\varphi_{gX}(Y)=d_Y\varphi_Y(g)\varphi_Y(X)/d_X$, which is $b(g,Y)\varphi_X(Y)$. This is the first identity below, and applying $\sigma$ to it gives the second:
\begin{equation}\label{eq:invertible}
\varphi_{g\otimes X}=b(g,-)\,\varphi_X\qquad\text{and}\qquad\hat\sigma(g\otimes X)=g^\ell\otimes\hat\sigma(X).
\end{equation}

\begin{lemma}\label{lem:symmetric}
Let $j\colon\mathcal{E}\to\mathcal{N}$ be a Galois-modular extension of a symmetric category $\mathcal{E}$, with (additive) central charge $c$. Then:
\begin{enumerate}[label=(\alph*)]
\item $\mathcal{N}$ is weakly integral, and $u_j(\sigma)\in\mathcal{E}_{\mathrm{pt}}[2]$ for every $\sigma$;
\item $\sigma(d_Z)=b(u_j(\sigma),Z)\,d_Z$ for all $Z\in\mathcal{O}(\mathcal{N})$;
\item $u_j$ is a group homomorphism, with kernel $\operatorname{Gal}(\overline{\mathbb{Q}}/K_\mathcal{N})$, where $K_\mathcal{N}=\mathbb{Q}(d_Z:Z\in\mathcal{O}(\mathcal{N}))$ is totally real;
\item $2c\in\mathbb{Z}$ and $\theta\circ u_j=\kappa^{2c}$.
\end{enumerate}
\end{lemma}

\begin{proof}
For the proof of (a), note that $\dim(\mathcal{N})=\dim(\mathcal{E})\dim(C_\mathcal{N}(\mathcal{E}))$ is an integer, and $d_{\hat\sigma(X)}^2=\sigma(d_X)^2\dim(\mathcal{N})/\sigma(\dim(\mathcal{N}))$ \cite[(4)]{MR4834521}, so $d_{u(\sigma)}=1$. By \cite[Theorem 4.1.6]{MR4834521}, $u(\sigma)\in\mathcal{E}$, and $u(\sigma)$ is self-dual because $\hat\sigma$ commutes with duality.

For the proof of (b), note that $\sigma\circ\varphi_\mathbf{1}=\varphi_{u(\sigma)}$ and $\varphi_u(Z)=s_{Z,u}=b(u,Z)\,d_Z$.

For the proof of (c), \eqref{eq:invertible} gives $u(\sigma\tau)=\hat\sigma(u(\tau))=u(\tau)^\ell\otimes u(\sigma)=u(\tau)\otimes u(\sigma)$, since $\ell$ is odd and $u(\tau)$ has order at most $2$. By (b), and because an invertible object with $b(u,-)\equiv1$ is transparent in the nondegenerate category $\mathcal{N}$, we have $u(\sigma)=\mathbf{1}$ if and only if $\sigma$ fixes every $d_Z$. Also by (b), $\sigma(d_Z)=\pm d_Z$ is real for every $\sigma$, so $K_\mathcal{N}$ is totally real.

For the proof of (d), recall that twists in $\mathcal{E}$ are $\pm1$. By \cite[Theorem II(iii)]{MR3435813}, $\sigma^2(\gamma^{-1}\theta_X)=\gamma^{-1}\theta_{\hat\sigma(X)}$ for $X\in\mathcal{O}(\mathcal{N})$, and taking $X=\mathbf{1}$ gives $\theta_{u(\sigma_\ell)}=\gamma/\sigma^2(\gamma)=\gamma^{1-\ell^2}$. Hence $\eta:=\gamma^2$ satisfies $\eta^{\ell^2-1}=1$ for every $\ell$ prime to the order $n$ of $\eta$, which is finite. So $(\mathbb{Z}/n)^\times$ has exponent at most $2$, which forces $n\mid24$. Writing $\gamma=\omega e^{2\pi ic/24}$ with $\omega^3=1$, this gives $2c\in\mathbb{Z}$. Since $\ell$ is prime to $6$, $24\mid\ell^2-1$, $(\ell^2-1)/24\equiv(\ell^2-1)/8\pmod2$, and
\[
\theta_{u(\sigma_\ell)}=e^{-2\pi ic(\ell^2-1)/24}=(-1)^{2c(\ell^2-1)/24}=\kappa(\sigma_\ell)^{2c}. \qedhere
\]
\end{proof}

% ======================================================================
\section{The structure theorem}\label{sec:structure}
% ======================================================================

In this section $\mathcal{B}$ is an arbitrary braided fusion category with symmetric center $\mathcal{E}$. Let $\iota:\mathcal{B}\to\mathcal{M}$ be a nondegenerate extension with centralizer $\mathcal{C}$, let $\mathcal{E}_1\subseteq\mathcal{E}$ be a fusion subcategory, and let $j:\mathcal{E}_1\to\mathcal{N}$ be a nondegenerate extension of $\mathcal{E}_1$. Define $\iota\oplus j:\mathcal{B}\to(\mathcal{M}\boxtimes\mathcal{N})^0_A$ by $X\mapsto F(\iota X\boxtimes\mathbf{1})$, where
\[
A:=\bigoplus_{V\in\mathcal{O}(\mathcal{E}_1)}\iota(V)\boxtimes j(V)^\ast\in\mathcal{M}\boxtimes\mathcal{N}.
\]
Here $A$ is the algebra $L_{\mathcal{E}_1}$ of Subsection \ref{subsec:relative}, transported along $\iota$ and $j$.

\begin{proposition}\label{prop:baer}
$\iota\oplus j$ is a nondegenerate extension of $\mathcal{B}$, with centralizer $F(\mathcal{C}\boxtimes C_\mathcal{N}(\mathcal{E}_1))\simeq\mathcal{C}\boxtimes_{\mathcal{E}_1}C_\mathcal{N}(\mathcal{E}_1)$. If the centralizers of $\iota$ and $j$ are integral, so is that of $\iota\oplus j$. If $\mathcal{B}$ is pseudounitary and $\iota,j$ are Galois-modular, so is $\iota\oplus j$, and its target has (additive) central charge $c(\mathcal{M})+c(\mathcal{N})$.
\end{proposition}

\begin{proof}
By \cite[Corollary 3.30, Proposition 5.4]{MR3039775}, the category $(\mathcal{M}\boxtimes\mathcal{N})^0_A$ is nondegenerate and Witt equivalent to $\mathcal{M}\boxtimes\mathcal{N}$. Since $\mathcal{E}_1\subseteq\mathcal{E}$ is transparent in $\mathcal{B}$, $\iota(\mathcal{B})\boxtimes\mathbf{1}\subseteq C_{\mathcal{M}\boxtimes\mathcal{N}}(A)$. By adjunction \cite[Theorem 1.6]{MR1936496}, for $X,Y\in\mathcal{O}(\mathcal{B})$,
\[
\operatorname{Hom}(F(X\boxtimes\mathbf{1}),F(Y\boxtimes\mathbf{1}))=\bigoplus_{V\in\mathcal{O}(\mathcal{E}_1)}\operatorname{Hom}(X,Y\otimes V)\otimes\operatorname{Hom}(\mathbf{1},j(V)^\ast)=\operatorname{Hom}(X,Y).
\]
So $\iota\oplus j$ is braided and fully faithful. The double centralizer theorem implies $C_\mathcal{M}(\mathcal{C})=\iota(\mathcal{B})$, so the symmetric center of $\mathcal{C}$ is $\mathcal{C}\cap\iota(\mathcal{B})=\mathcal{E}$, and $\mathcal{E}_1$ is transparent in $\mathcal{C}$. Also $j(\mathcal{E}_1)$ is transparent in $C_\mathcal{N}(\mathcal{E}_1)$. Hence $\mathcal{C}\boxtimes C_\mathcal{N}(\mathcal{E}_1)$ contains $A$, lies in $C_{\mathcal{M}\boxtimes\mathcal{N}}(A)$, and centralizes $\iota(\mathcal{B})\boxtimes\mathbf{1}$. Its image under $F$ is $\mathcal{C}\boxtimes_{\mathcal{E}_1}C_\mathcal{N}(\mathcal{E}_1)$, with Frobenius--Perron dimension
\[
\mathrm{FPdim}(\mathcal{C})\mathrm{FPdim}(\mathcal{N})/\mathrm{FPdim}(\mathcal{E}_1)^2=\mathrm{FPdim}((\mathcal{M}\boxtimes\mathcal{N})^0_A)/\mathrm{FPdim}(\mathcal{B}).
\]
This is the Frobenius--Perron dimension of the centralizer of $(\iota\oplus j)(\mathcal{B})$, so the two coincide. For modules $L$ over a connected \'etale algebra $A$ in a braided fusion category $\mathcal{P}$, the functions $L\mapsto\mathrm{FPdim}(L)$ and $L\mapsto\mathrm{FPdim}_\mathcal{P}(L)/\mathrm{FPdim}(A)$ are positive eigenvectors of the action $X\cdot L=X\otimes L$ of $K(\mathcal{P})$ on the indecomposable module category $\mathcal{P}_A$. So they are proportional, and they agree at $L=A$. If $\mathcal{C}$ and $C_\mathcal{N}(\mathcal{E}_1)$ are integral, applying this to $A\in\mathcal{C}\boxtimes C_\mathcal{N}(\mathcal{E}_1)$ shows that the simple objects of the centralizer have Frobenius--Perron dimension equal to an integer divided by $\mathrm{FPdim}(\mathcal{E}_1)$. These are rational algebraic integers, hence integers.

Finally, let $\mathcal{B}$ be pseudounitary and let $\iota$ and $j$ be Galois-modular. Then the centralizer of $\iota\oplus j$ is integral, hence pseudounitary \cite[Proposition 8.24]{MR2183279}, so the target of $\iota\oplus j$ is pseudounitary \cite[Proposition 2.4]{ART_2026__3_3_579_0}. Witt equivalent pseudounitary categories have the same central charge \cite[Lemma 5.27]{MR3039775}, so this central charge is $c(\mathcal{M})+c(\mathcal{N})$.
\end{proof}

\begin{note}\label{rem:baer}
The name \emph{Baer sum} comes from Baer sums of group extensions, and for $\mathcal{E}_1=\mathcal{E}\simeq\mathrm{Rep}(G)$ the analogy is literal. By de-equivariantization, an element of $\mathrm{Mext}(\mathrm{Rep}(G))$ is the same as a $G$-graded pointed category $\mathrm{Vec}_G^\omega$, that is, a categorical group $\mathcal{G}_\omega$ with $\pi_0=G$, $\pi_1=\mathbb{C}^\times$ and $k$-invariant $\omega\in H^3(G,\mathbb{C}^\times)$ \cite{MR2068521,MR3242743}. Here $L_\mathcal{E}\cong\mathrm{Fun}((G\times G)/\Delta)$, so condensing it restricts the symmetry of $\mathcal{Z}(\mathrm{Vec}_G^\omega)\boxtimes\mathcal{Z}(\mathrm{Vec}_G^{\omega'})$ from $G\times G$ to the diagonal. Thus $\boxtimes_\mathcal{E}$ pulls $\mathcal{G}_\omega\times\mathcal{G}_{\omega'}$ back along $G\to G\times G$ and identifies the two copies of $\mathbb{C}^\times$ by multiplication. This is the Baer sum of $2$-group extensions, and its $k$-invariant is $\omega\omega'$. More generally, let $\mathcal{D}$ be nondegenerate with a braided $G$-action. The minimal nondegenerate extensions of $\mathcal{D}^G$ form an $H^3(G,\mathbb{C}^\times)$-torsor once an obstruction in $H^4(G,\mathbb{C}^\times)$ vanishes \cite{MR2677836}, and the same computation shows that $\omega$ acts by Baer sums with $\mathcal{Z}(\mathrm{Vec}_G^\omega)$. This is the theory of group extensions with non-abelian kernel $N$ one degree up: there the obstruction lies in $H^3(G,Z(N))$, and $H^2(G,Z(N))$ acts by Baer sums with central extensions \cite{MR20996}.
\end{note}

We can now state and prove the structure theorem.

\vspace{5 mm}

\begin{theorem}\label{thm:structure}
Let $\mathcal{B}$ be a braided fusion category with symmetric center $\mathcal{E}$.
\begin{enumerate}[label=(\alph*)]
\item Let $\iota:\mathcal{B}\to\mathcal{M}$ be a nondegenerate extension with centralizer $\mathcal{C}$. Then $\mathcal{E}$ is also the symmetric center of $\mathcal{C}$, $\iota(\mathcal{B})\vee\mathcal{C}\simeq\mathcal{B}\boxtimes_\mathcal{E}\mathcal{C}$, and $\mathcal{M}$ is a minimal nondegenerate extension of $\mathcal{B}\boxtimes_\mathcal{E}\mathcal{C}$. Conversely, let $\mathcal{C}$ be a braided fusion category with symmetric center $\mathcal{E}$. Then every nondegenerate $\mathcal{M}\supseteq\mathcal{B}\boxtimes_{\mathcal{E}}\mathcal{C}$ with $C_\mathcal{M}(\mathcal{B}\boxtimes_{\mathcal{E}}\mathcal{C})=\mathcal{E}$ is a nondegenerate extension of $\mathcal{B}$ with centralizer $\mathcal{C}$.
\item Any two pseudounitary nondegenerate extensions $\iota,\iota'$ of $\mathcal{B}$ are stably equivalent: there are nondegenerate extensions $j,j'$ of $\mathcal{E}$ with $\iota\oplus j\simeq\iota'\oplus j'$. If $\iota$ and $\iota'$ have integral centralizers, $j$ and $j'$ can be chosen to have integral centralizers.
\item If $\mathcal{B}$ is pseudounitary with a minimal nondegenerate extension $\iota_0$, then for every pseudounitary braided fusion category $\mathcal{C}$ with symmetric center identified with $\mathcal{E}$, $\mathcal{N}\mapsto\iota_0\oplus(\mathcal{E}\subseteq\mathcal{N})$ is a bijection from $\mathrm{Mext}(\mathcal{C})$ to $\mathrm{Mext}(\mathcal{B}\boxtimes_\mathcal{E}\mathcal{C})$, i.e.\ to the nondegenerate extensions of $\mathcal{B}$ with centralizer identified with $\mathcal{C}$, up to equivalences respecting $\mathcal{B}\boxtimes_\mathcal{E}\mathcal{C}$.
\end{enumerate}
\end{theorem}

In (c), and whenever we consider extensions with centralizer $\mathcal{C}$, the identification of $C_\mathcal{M}(\iota\mathcal{B})$ with $\mathcal{C}$ from (a) is part of the data; forgetting it can identify extensions. Part (c) says that, for pseudounitary $\mathcal{B},\mathcal{C}$ and when a minimal extension exists, every nondegenerate extension of $\mathcal{B}$ with centralizer $\mathcal{C}$ is obtained by gluing it along $\mathcal{E}$ to a minimal extension of $\mathcal{C}$. Part (b) needs no existence hypothesis. We use freely that Baer sums along $\mathcal{E}$ are associative, commutative and functorial up to equivalence of extensions, and that $\iota\oplus(\mathcal{E}\subseteq\mathcal{Z}(\mathcal{E}))\simeq\iota$. By (a), every nondegenerate extension of $\mathcal{B}$ is a minimal nondegenerate extension of some $\mathcal{B}\boxtimes_\mathcal{E}\mathcal{C}$ (cf.\ the proof of \cite[Proposition 2.1]{ART_2026__3_3_579_0}), and for these the statements are \cite[Propositions 4.13 and 4.14, Lemma 4.18]{MR3613518}.

\begin{proof}[Proof of Theorem~\ref{thm:structure}]
For the proof of (a), by the double centralizer theorem, $C_\mathcal{M}(\mathcal{C})=\iota(\mathcal{B})$. So $\mathcal{C}'=\mathcal{C}\cap\iota(\mathcal{B})=\iota(\mathcal{B})'=\mathcal{E}$. Lemma~\ref{lem:join} then gives $\iota(\mathcal{B})\vee\mathcal{C}\simeq\mathcal{B}\boxtimes_{\mathcal{E}}\mathcal{C}$, and moreover $C_\mathcal{M}(\iota\mathcal{B}\vee\mathcal{C})=C_\mathcal{M}(\iota\mathcal{B})\cap C_\mathcal{M}(\mathcal{C})=\mathcal{C}\cap\iota(\mathcal{B})=\mathcal{E}$.

Conversely, let $\mathcal{D}:=\mathcal{B}\boxtimes_{\mathcal{E}}\mathcal{C}\subseteq\mathcal{M}$ with $C_\mathcal{M}(\mathcal{D})=\mathcal{E}$. Then $C_\mathcal{M}(\mathcal{B})\supseteq\mathcal{C}$, and
\[
\mathrm{FPdim}(C_\mathcal{M}(\mathcal{B}))=\frac{\mathrm{FPdim}(\mathcal{M})}{\mathrm{FPdim}(\mathcal{B})}=\frac{\mathrm{FPdim}(\mathcal{D})\mathrm{FPdim}(\mathcal{E})}{\mathrm{FPdim}(\mathcal{B})}=\mathrm{FPdim}(\mathcal{C}).
\]

Now for the proof of (b), let $\mathcal{C}_1,\mathcal{C}_2$ be the centralizers of $\iota\colon\mathcal{B}\to\mathcal{M}_1$ and $\iota'\colon\mathcal{B}\to\mathcal{M}_2$.

\emph{Step 1:} The canonical algebra $L_\mathcal{B}=\bigoplus_{X\in\mathcal{O}(\mathcal{B})}X^*\boxtimes X$ is a connected \'etale algebra in $\mathcal{B}^{\mathrm{rev}}\boxtimes\mathcal{B}\subseteq\mathcal{M}_1^{\mathrm{rev}}\boxtimes\mathcal{M}_2$ \cite[Lemma 3.5]{MR3039775}. Let $\mathcal{N}:=(\mathcal{M}_1^{\mathrm{rev}}\boxtimes\mathcal{M}_2)^0_{L_\mathcal{B}}$; it is nondegenerate of Frobenius--Perron dimension $\mathrm{FPdim}(\mathcal{C}_1)\mathrm{FPdim}(\mathcal{C}_2)$ \cite[Corollaries 3.30 and 3.32]{MR3039775}, and pseudounitary (Lemma~\ref{lem:trace}).

Objects of $\mathcal{E}\boxtimes\mathbf{1}$ and of $\mathcal{C}_1^{\mathrm{rev}}\boxtimes\mathcal{C}_2$ centralize $L_\mathcal{B}$. As in the proof of Proposition~\ref{prop:baer}, $V\mapsto F(V\boxtimes\mathbf{1})$ embeds $\mathcal{E}$ in $\mathcal{N}$. The summands of $L_\mathcal{B}$ in $\mathcal{E}\boxtimes\mathcal{E}$ form a subalgebra isomorphic to $L_\mathcal{E}$, and induction along $L_\mathcal{E}\subseteq L_\mathcal{B}$ is a braided tensor functor $\mathcal{C}_1^{\mathrm{rev}}\boxtimes_{\mathcal{E}}\mathcal{C}_2\to\mathcal{N}$ with image $F(\mathcal{C}_1^{\mathrm{rev}}\boxtimes\mathcal{C}_2)$. Only these summands contribute to Hom spaces between objects of $\mathcal{C}_1^{\mathrm{rev}}\boxtimes\mathcal{C}_2$, so this functor is fully faithful. Its image has Frobenius--Perron dimension $\mathrm{FPdim}(\mathcal{C}_1)\mathrm{FPdim}(\mathcal{C}_2)/\mathrm{FPdim}(\mathcal{E})=\mathrm{FPdim}(C_\mathcal{N}(\mathcal{E}))$ and centralizes $\mathcal{E}$. Hence $C_\mathcal{N}(\mathcal{E})\simeq\mathcal{C}_1^{\mathrm{rev}}\boxtimes_{\mathcal{E}}\mathcal{C}_2$. So $j_\mathcal{N}\colon\mathcal{E}\subseteq\mathcal{N}$ is a nondegenerate extension of $\mathcal{E}$, with integral centralizer if $\mathcal{C}_1$ and $\mathcal{C}_2$ are integral (as in the proof of Proposition~\ref{prop:baer}).

\emph{Step 2:} In $\mathcal{Z}(\mathcal{C}_1)$ the images of $\mathcal{C}_1$ and $\mathcal{C}_1^{\mathrm{rev}}$ centralize each other and meet in $\mathcal{C}_1'=\mathcal{E}$ \cite[Proposition 3.7, Remark 3.8]{MR2609644}. By Lemma~\ref{lem:join} their join is $\mathcal{C}_1\boxtimes_{\mathcal{E}}\mathcal{C}_1^{\mathrm{rev}}$, of Frobenius--Perron dimension $\mathrm{FPdim}(\mathcal{Z}(\mathcal{C}_1))/\mathrm{FPdim}(\mathcal{E})$, so $C_{\mathcal{Z}(\mathcal{C}_1)}(\mathcal{E})\simeq\mathcal{C}_1\boxtimes_{\mathcal{E}}\mathcal{C}_1^{\mathrm{rev}}$. Thus $j_{\mathcal{Z}(\mathcal{C}_1)}\colon\mathcal{E}\subseteq\mathcal{Z}(\mathcal{C}_1)$ is a nondegenerate extension of $\mathcal{E}$, with integral centralizer if $\mathcal{C}_1$ is integral.

\emph{Step 3:} By Proposition~\ref{prop:baer}, the centralizers of $\mathcal{B}$ in $\iota\oplus j_\mathcal{N}$ and in $\iota'\oplus j_{\mathcal{Z}(\mathcal{C}_1)}$ are $\mathcal{C}_1\boxtimes_{\mathcal{E}}(\mathcal{C}_1^{\mathrm{rev}}\boxtimes_{\mathcal{E}}\mathcal{C}_2)$ and $\mathcal{C}_2\boxtimes_{\mathcal{E}}(\mathcal{C}_1\boxtimes_{\mathcal{E}}\mathcal{C}_1^{\mathrm{rev}})$. Both are equivalent to $\widetilde{\mathcal{C}}:=\mathcal{C}_1\boxtimes_{\mathcal{E}}\mathcal{C}_1^{\mathrm{rev}}\boxtimes_{\mathcal{E}}\mathcal{C}_2$ via associativity and commutativity of $\boxtimes_\mathcal{E}$. These respect the embeddings of $\mathcal{E}$, since in a relative Deligne product the two copies of $\mathcal{E}$ are identified: $F(V\boxtimes\mathbf{1})\cong F(\mathbf{1}\boxtimes V)$ for $V\in\mathcal{O}(\mathcal{E})$. By (a), both are minimal nondegenerate extensions of the pseudounitary category $\mathcal{B}\boxtimes_{\mathcal{E}}\widetilde{\mathcal{C}}$, and they restrict to the given embedding of $\mathcal{B}$. By Theorem~\ref{thm:LKW}(b) there is $\mathcal{K}\in\mathrm{Mext}(\mathcal{E})$ with $\iota\oplus j_\mathcal{N}\simeq(\iota'\oplus j_{\mathcal{Z}(\mathcal{C}_1)})\oplus j_\mathcal{K}$. By associativity this equals $\iota'\oplus j'$ with $j'=j_{\mathcal{Z}(\mathcal{C}_1)\boxtimes_{\mathcal{E}}\mathcal{K}}$. Every element of $\mathrm{Mext}(\mathcal{E})$ has integral centralizer $\mathcal{E}$.

Lastly, for the proof of (c), let $\iota$ have centralizer $\mathcal{C}$, and apply Step~1 to $(\iota_0,\iota)$. This gives $\mathcal{N}$ with $C_\mathcal{N}(\mathcal{E})\simeq\mathcal{C}$, so $C_\mathcal{N}(\mathcal{C})=\mathcal{E}$ and $\mathcal{N}\in\mathrm{Mext}(\mathcal{C})$. By Proposition~\ref{prop:baer}, $\iota_0\oplus j_\mathcal{N}$ has centralizer $\mathcal{E}\boxtimes_{\mathcal{E}}\mathcal{C}\simeq\mathcal{C}$, as does $\iota$. By (a) and Theorem~\ref{thm:LKW}(b), $\iota\simeq(\iota_0\oplus j_\mathcal{N})\oplus j_\mathcal{K}=\iota_0\oplus j_{\mathcal{N}\boxtimes_{\mathcal{E}}\mathcal{K}}$ for some $\mathcal{K}\in\mathrm{Mext}(\mathcal{E})$, and $\mathcal{N}\boxtimes_{\mathcal{E}}\mathcal{K}\in\mathrm{Mext}(\mathcal{C})$. So the map is surjective. Both sides are empty or $\mathrm{Mext}(\mathcal{E})$-torsors by Theorem~\ref{thm:LKW}(b), and the map is equivariant by associativity, so it is bijective.
\end{proof}

% ======================================================================
\section{Galois actions}\label{sec:galois}
% ======================================================================

\begin{proposition}\label{prop:formula}
Let $\iota$ be Galois-modular, and let $j:\mathcal{E}_1\to\mathcal{N}$ be a Galois-modular extension of a fusion subcategory $\mathcal{E}_1\subseteq\mathcal{E}$. Then for every $\sigma\in\operatorname{Gal}(\overline{\mathbb{Q}}/\mathbb{Q})$ and $X\in\mathcal{O}(\mathcal{B})$,
\[
\hat\sigma_{\iota\oplus j}(X)=\hat\sigma_\iota(X)\otimes u_j(\sigma).
\]
\end{proposition}

\begin{proof}
Let $u=u_j(\sigma)$. By Lemma~\ref{lem:symmetric}, $u\in(\mathcal{E}_1)_{\mathrm{pt}}[2]\subseteq\mathcal{E}$, and $\sigma(d_Z)=b(u,Z)\,d_Z$ for $Z\in\mathcal{O}(\mathcal{N})$. Let $L$ be a simple local $A$-module with underlying object $\bigoplus_iY_i\boxtimes Z_i$. By Lemma~\ref{lem:trace},
\[
\varphi_{F(X\boxtimes\mathbf{1})}(L)=\frac{s_{L,F(X\boxtimes\mathbf{1})}}{d_X}=\frac1{\dim(\mathcal{E}_1)}\sum_i\varphi_X(Y_i)\,d_{Z_i}.
\]
Applying $\sigma$ replaces $\varphi_X(Y_i)$ by $\varphi_{\hat\sigma_\iota X}(Y_i)$ and $d_{Z_i}$ by $b(u,Z_i)\,d_{Z_i}$. The object $u\boxtimes j(u)^{-1}$ is an invertible summand of $A$, so Lemma~\ref{lem:trace} gives $b(u,Y_i)=b(u,Z_i)$. With \eqref{eq:invertible} this yields
\[
\sigma\bigl(\varphi_{F(X\boxtimes\mathbf{1})}(L)\bigr)=\frac1{\dim(\mathcal{E}_1)}\sum_i\varphi_{\hat\sigma_\iota X\otimes u}(Y_i)\,d_{Z_i}=\varphi_{F((\hat\sigma_\iota X\otimes u)\boxtimes\mathbf{1})}(L).
\]
As $L$ is arbitrary, $\sigma\circ\varphi_{F(X\boxtimes\mathbf{1})}=\varphi_{F((\hat\sigma_\iota X\otimes u)\boxtimes\mathbf{1})}$.
\end{proof}

For $g\in\mathcal{E}_{\mathrm{pt}}[2]$ and a character $\chi\colon\operatorname{Gal}(\overline{\mathbb{Q}}/\mathbb{Q})\to\{\pm1\}$, write $g^\chi$ for the map sending $\sigma$ to $g$ if $\chi(\sigma)=-1$ and to $\mathbf{1}$ otherwise.

\begin{example}\label{ex:generators}
\leavevmode
\begin{enumerate}[label=(\roman*)]
\item For an odd prime $p$, the metaplectic category $SO(p)_2$ has simple objects $\mathbf{1},Z$ of dimension $1$, $(p-1)/2$ objects of dimension $2$, and two of dimension $\sqrt p$ \cite{MR3541678}. Its integral simple objects form the adjoint subcategory, which is the centralizer of $\{\mathbf{1},Z\}\simeq\mathrm{Rep}(\mathbb{Z}/2)$ \cite[Corollary 6.9]{MR2383894}. By Lemma~\ref{lem:symmetric}(b), $u(\sigma)=Z$ if and only if $\sigma(\sqrt p)=-\sqrt p$.
\item Write $\mathbf{1},\mu,\psi$ for the simple objects of $\mathrm{Ising}$, with $d_\mu=\sqrt2$. Then $\{\mathbf{1},\psi\boxtimes\psi\}\simeq\mathrm{Rep}(\mathbb{Z}/2)$ in $\mathrm{Ising}^{\boxtimes2}$ has integral centralizer, with simple objects $\mathbf{1}$, $\psi\boxtimes\mathbf{1}$, $\mathbf{1}\boxtimes\psi$, $\psi\boxtimes\psi$ and $\mu\boxtimes\mu$, and $u=(\psi\boxtimes\psi)^\kappa$.
\item $\{\mathbf{1},\psi\}\simeq\mathrm{sVec}$ in $\mathrm{Ising}$ has integral centralizer $\{\mathbf{1},\psi\}$, and $u=\psi^\kappa$.
\end{enumerate}
\end{example}

We now restate the main results of Section~\ref{sec:tech} and prove them.

\begin{theorem}\label{thm:main}
Let $\mathcal{B}$ be a pseudounitary braided fusion category with symmetric center $\mathcal{E}$.
\begin{enumerate}[label=(\alph*)]
\item If $\iota,j$ are Galois-modular extensions of $\mathcal{B}$ and a fusion subcategory $\mathcal{E}_1\subseteq\mathcal{E}$, respectively, then $\iota\oplus j$ is Galois-modular, and $\hat\sigma_{\iota\oplus j}(X)=\hat\sigma_\iota(X)\otimes u_j(\sigma)$ for all $X\in\mathcal{O}(\mathcal{B})$.
\item Under Baer sums along $\mathcal{E}$, the integral equivalence classes of Galois-modular extensions of $\mathcal{E}$ form an abelian group. The map $j\mapsto u_j$ is an isomorphism from it onto the set $\mathcal{U}(\mathcal{E})$ of all $u_j$, a group under pointwise product. The integral equivalence classes of Galois-modular extensions of $\mathcal{B}$ form a torsor over this group, and $\iota\mapsto\hat\sigma_\iota$ is an equivariant bijection onto $\mathcal{A}(\mathcal{B})$. In particular, $\mathcal{A}(\mathcal{B})$ is a $\mathcal{U}(\mathcal{E})$-torsor under $(\hat\sigma_\iota,u)\mapsto\hat\sigma_\iota\otimes u$, and two Galois-modular extensions of $\mathcal{B}$ induce the same action if and only if they are integrally equivalent.
\item $\mathcal{U}(\mathcal{E})$ is the group of homomorphisms $u\colon\operatorname{Gal}(\overline{\mathbb{Q}}/\mathbb{Q})\to\mathcal{E}_{\mathrm{pt}}[2]$ that factor through $\operatorname{Gal}(K/\mathbb{Q})$ for some real multiquadratic field $K$ and satisfy $\theta\circ u\in\{1,\kappa\}$, where $\theta$ is the twist.
\end{enumerate}
\end{theorem}

\begin{proof}
Part (a) is Propositions~\ref{prop:baer} and~\ref{prop:formula}.

For the proof of (b), recall from Lemma~\ref{lem:symmetric}(c) that a Galois-modular extension $z$ of $\mathcal{E}$ is integral if and only if $u_z\equiv\mathbf{1}$. Baer sums are associative and commutative, and by (a), $u_{j\oplus j'}=u_ju_{j'}$ and $\hat\sigma_{\iota\oplus j}=\hat\sigma_\iota\otimes u_j$. Hence integral equivalence is an equivalence relation compatible with Baer sums, and $j\mapsto u_j$ and $\iota\mapsto\hat\sigma_\iota$ are constant on classes. So the classes of Galois-modular extensions of $\mathcal{E}$ form a commutative monoid $\mathcal I$, whose identity is the class of $\mathcal{E}\subseteq\mathcal{Z}(\mathcal{E})$, since Baer sum with it is trivial (see the paragraph after Theorem~\ref{thm:structure}). The map $u\colon\mathcal I\to\mathcal{U}(\mathcal{E})$ is a surjective homomorphism, and $\mathcal I$ acts on the classes of extensions of $\mathcal{B}$, compatibly with $\hat\sigma$ and $u$.

Next let $\iota,\iota'$ be Galois-modular, with centralizer $\mathcal{C}_1$ of $\iota$. By Step~3 of the proof of Theorem~\ref{thm:structure}(b), $\iota\oplus j_\mathcal{N}\simeq\iota'\oplus j_{\mathcal{Z}(\mathcal{C}_1)}\oplus j_\mathcal{K}$ with $\mathcal{K}\in\mathrm{Mext}(\mathcal{E})$, and these extensions of $\mathcal{E}$ are Galois-modular; for $j_\mathcal{N}$ this holds because its centralizer is integral, by \cite[Proposition 8.24]{MR2183279} and \cite[Proposition 2.4]{ART_2026__3_3_579_0}. Since $\mathcal{Z}(\mathcal{C}_1)$ is integral, the classes satisfy $[\iota]+[j_\mathcal{N}]=[\iota']+[j_\mathcal{K}]$, and (a) gives $\hat\sigma_{\iota'}=\hat\sigma_\iota\otimes u_{j_\mathcal{N}}u_{j_\mathcal{K}}$.

The action of $\mathcal{U}(\mathcal{E})$ on $\mathcal{A}(\mathcal{B})$ is free: if $\hat\sigma_\iota(\mathbf{1})\otimes u(\sigma)\cong\hat\sigma_\iota(\mathbf{1})$, apply $\widehat{\sigma^{-1}}$ and \eqref{eq:invertible}: this gives $u(\sigma)^{\ell^{-1}}\cong\mathbf{1}$, so $u(\sigma)=\mathbf{1}$. With the previous paragraph, $\mathcal{U}(\mathcal{E})$ acts freely and transitively on $\mathcal{A}(\mathcal{B})$ by $\hat\sigma_\iota\otimes u_j=\hat\sigma_{\iota\oplus j}$, and $\mathcal{A}(\mathcal{B})\ne\emptyset$ by \cite{ART_2026__3_3_579_0}.

The map $u\colon\mathcal I\to\mathcal{U}(\mathcal{E})$ is injective. Let $j,j'$ be Galois-modular extensions of $\mathcal{E}$ with $u_j=u_{j'}$, and let $\mathcal{C}_1$ be the centralizer of $j\colon\mathcal{E}\to\mathcal{N}$. Apply the proof of Theorem~\ref{thm:structure}(b) to $(\mathcal{E},j,j')$. Step~1 gives the Baer sum $j^{\mathrm{rev}}\oplus j'$, and Step~3 gives $j\oplus(j^{\mathrm{rev}}\oplus j')\simeq j'\oplus j_{\mathcal{Z}(\mathcal{C}_1)}\oplus j_\mathcal{K}$ for $\mathcal{K}\in\mathrm{Mext}(\mathcal{E})$. By Lemma~\ref{lem:symmetric}(b), $b(u_{j^{\mathrm{rev}}}(\sigma),Z)=b(u_j(\sigma),Z)$ for all $Z\in\mathcal{O}(\mathcal{N})$, since $\mathcal{N}^{\mathrm{rev}}$ has the same dimensions as $\mathcal{N}$ and $b$ is real on $\mathcal{E}_{\mathrm{pt}}[2]$. As $\mathcal{N}$ is nondegenerate, $u_{j^{\mathrm{rev}}}=u_j$. As every $u$ has order at most $2$, comparing both sides gives $u_{j_\mathcal{K}}\equiv\mathbf{1}$, and $u_{j^{\mathrm{rev}}\oplus j'}=u_ju_{j'}\equiv\mathbf{1}$. So the extensions $j^{\mathrm{rev}}\oplus j'$, $j_{\mathcal{Z}(\mathcal{C}_1)}$ and $j_\mathcal{K}$ are integral, and $j,j'$ are integrally equivalent.

Thus $u\colon\mathcal I\to\mathcal{U}(\mathcal{E})$ is an isomorphism, and $\mathcal I$ is a group. By the relation $[\iota]+[j_\mathcal{N}]=[\iota']+[j_\mathcal{K}]$ above, $[\iota']=[\iota]+[j_\mathcal{N}]-[j_\mathcal{K}]$, so $\mathcal I$ acts transitively on the classes of extensions of $\mathcal{B}$. If $\hat\sigma_\iota=\hat\sigma_{\iota'}$, then $u_{j_\mathcal{N}}u_{j_\mathcal{K}}\equiv\mathbf{1}$ by freeness, so $[j_\mathcal{N}]=[j_\mathcal{K}]$ and $[\iota]=[\iota']$. Hence $\iota\mapsto\hat\sigma_\iota$ is injective on classes, and it is surjective by definition. Since it is equivariant along the isomorphism $u$ and $\mathcal{A}(\mathcal{B})$ is a $\mathcal{U}(\mathcal{E})$-torsor, the classes of extensions of $\mathcal{B}$ form an $\mathcal I$-torsor.

For the proof of (c), let $\mathcal{G}$ be the group described there. By Lemma~\ref{lem:symmetric}, every $u_j$ factors through $\operatorname{Gal}(K_\mathcal{N}/\mathbb{Q})$ with $K_\mathcal{N}$ totally real. This group embeds in $\mathcal{E}_{\mathrm{pt}}[2]$, so it is an elementary abelian $2$-group, $K_\mathcal{N}$ is multiquadratic, and $u_j\in\mathcal{G}$.

Conversely, let $u\in\mathcal{G}$. Since $\mathcal{E}$ is symmetric, $\theta$ restricts to a homomorphism $\mathcal{E}_{\mathrm{pt}}[2]\to\{\pm1\}$. Choose an $\mathbb{F}_2$-basis $g_1,\dots,g_r$ of its kernel, together with one $f$ with $\theta_f=-1$ if $\theta\ne1$. Write $u=\prod_ig_i^{a_i}\,f^{\chi_f}$ with characters $a_i,\chi_f$ of order at most $2$; then $\chi_f=\theta\circ u\in\{1,\kappa\}$. Each $a_i$ factors through $\operatorname{Gal}(K/\mathbb{Q})$, so it is trivial or the character of a quadratic subfield $\mathbb{Q}(\sqrt m)$ of $K$, which is real. That character is the product of the characters of $\mathbb{Q}(\sqrt p)$ over the primes $p\mid m$.

Start with $\mathcal{E}\subseteq\mathcal{Z}(\mathcal{E})$, which has $u\equiv\mathbf{1}$ since $\mathcal{Z}(\mathcal{E})$ is integral. Baer-sum it successively along the subcategories $\{\mathbf{1},g_i\}\simeq\mathrm{Rep}(\mathbb{Z}/2)$ with the extensions of Example~\ref{ex:generators}(i),(ii), and, if $\chi_f=\kappa$, along $\{\mathbf{1},f\}\simeq\mathrm{sVec}$ with Example~\ref{ex:generators}(iii). By (a), the result is a Galois-modular extension $j$ of $\mathcal{E}$ with $u_j=u$. Hence $\mathcal{U}(\mathcal{E})=\mathcal{G}$.
\end{proof}

\begin{corollary}\label{cor:A}
Let $G$ be a finite group. Every induced Galois action on $\mathrm{Rep}(G)$ has the form $X\mapsto u(\sigma)\otimes(\sigma\circ X)$. Here $\sigma\circ X$ is the Galois conjugate representation, and $u\colon\operatorname{Gal}(\overline{\mathbb{Q}}/\mathbb{Q})\to\operatorname{Hom}(G,\{\pm1\})$ is a homomorphism that factors through $\operatorname{Gal}(K/\mathbb{Q})$ for some real multiquadratic field $K$. Every such $u$ occurs.
\end{corollary}

\begin{proof}
In $\mathcal{Z}(\mathrm{Vec}_G)$, the $(e,\psi)$ form $\mathrm{Rep}(G)$, with integral centralizer. By the Coste--Gannon--Ruelle formula \cite{MR1770077}, $\varphi_{(e,\psi)}(b,\beta)=|\mathrm{cl}(b)|\,\dim\beta\;\overline{\psi(b)}/\psi(1)$. Write $\hat\sigma(e,\psi)=(e,\psi')$; then $\psi'(1)=\psi(1)$ by \cite[(4)]{MR4834521}, since $\mathcal{Z}(\mathrm{Vec}_G)$ is integral. Comparing values for $(b,\mathbf{1})$ gives $\psi'=\sigma\circ\psi$. So this extension induces Galois conjugation. Now apply parts (b) and (c) of Theorem~\ref{thm:main}, noting that $\mathrm{Rep}(G)_{\mathrm{pt}}[2]=\operatorname{Hom}(G,\{\pm1\})$ and that $\theta=1$ on $\mathrm{Rep}(G)$.
\end{proof}

\begin{corollary}\label{cor:B}
Let $\mathcal{B}$ be a pseudounitary braided fusion category with symmetric center $\mathcal{E}$. If $\iota,\iota'$ are Galois-modular extensions of $\mathcal{B}$ with (additive) central charges $c,c'$, then $2(c-c')\in\mathbb{Z}$, and $c-c'\in\mathbb{Z}$ if $\mathcal{E}$ is Tannakian. If $\mathcal{E}$ is Tannakian and $\mathcal{B}$ has minimal modular extensions, then all minimal modular extensions induce the same Galois action on $\mathcal{B}$.
\end{corollary}

\begin{proof}
In the proof of Theorem~\ref{thm:structure}(b), $\mathcal{N}$ is Witt equivalent to $\mathcal{M}_1^{\mathrm{rev}}\boxtimes\mathcal{M}_2$ \cite[Proposition 5.4]{MR3039775}, so $c(\mathcal{N})=c'-c$ \cite[Lemma 5.27]{MR3039775}. By Lemma~\ref{lem:symmetric}(d), $2c(\mathcal{N})\in\mathbb{Z}$. If $\mathcal{E}$ is Tannakian, then $\theta\circ u_{j_\mathcal{N}}=1$, so $c(\mathcal{N})\in\mathbb{Z}$. Now let $\mathcal{E}\simeq\mathrm{Rep}(G)$ be Tannakian. Two minimal extensions differ by some $\mathcal{K}=\mathcal{Z}(\mathrm{Vec}_G^\omega)$ (Theorem~\ref{thm:LKW}), which is integral, so $u_{j_\mathcal{K}}\equiv\mathbf{1}$. Apply part (a) of Theorem~\ref{thm:main}.
\end{proof}

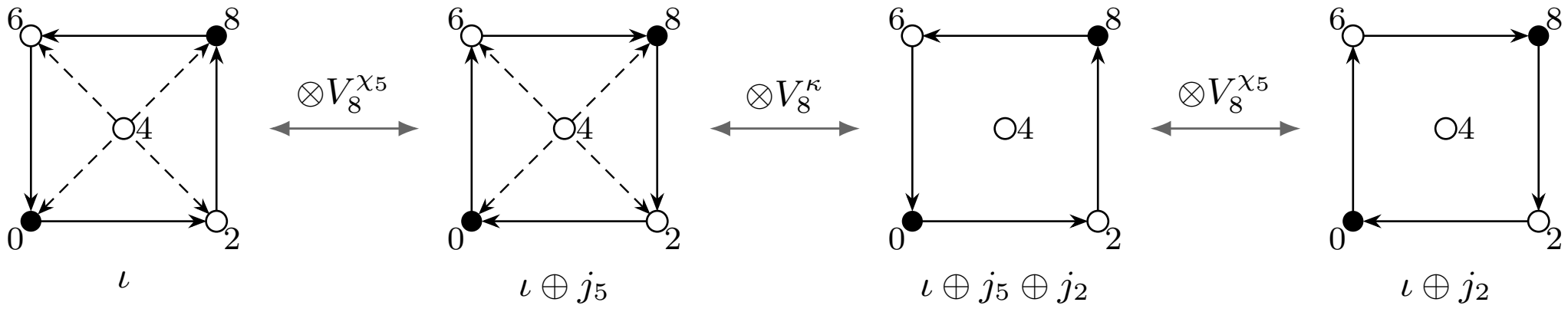
\begin{figure}[ht]
\centering
\begin{tikzpicture}[
  s3/.style={-{Stealth[length=1.6mm,width=1.2mm]}, line width=.5pt, shorten >=2.4pt, shorten <=2.4pt},
  s11/.style={{Stealth[length=1.4mm,width=1.1mm]}-{Stealth[length=1.4mm,width=1.1mm]}, densely dashed, line width=.45pt, shorten >=2.4pt, shorten <=2.4pt},
  eobj/.style={circle, fill=black, inner sep=0pt, minimum size=1.8mm},
  obj/.style={circle, fill=white, draw=black, line width=.5pt, inner sep=0pt, minimum size=1.8mm},
  lab/.style={font=\scriptsize, inner sep=1.5pt},
  cay/.style={{Latex[length=1.8mm]}-{Latex[length=1.8mm]}, line width=.55pt, black!60}]
\newcommand{\tpanel}[5]{%
 \begin{scope}[shift={(#1,#2)}]
  \coordinate (p0) at (0,0); \coordinate (p2) at (1.6,0); \coordinate (p8) at (1.6,1.6);
  \coordinate (p6) at (0,1.6); \coordinate (p4) at (0.8,0.8);
  \ifnum#3=1 \draw[s3] (p0)--(p2); \draw[s3] (p2)--(p8); \draw[s3] (p8)--(p6); \draw[s3] (p6)--(p0);
  \else \draw[s3] (p2)--(p0); \draw[s3] (p8)--(p2); \draw[s3] (p6)--(p8); \draw[s3] (p0)--(p6); \fi
  \ifnum#4=1 \draw[s11] (p0)--(p8); \draw[s11] (p2)--(p6); \fi
  \node[eobj] at (p0) {}; \node[eobj] at (p8) {};
  \node[obj] at (p2) {}; \node[obj] at (p6) {}; \node[obj] at (p4) {};
  \node[lab, anchor=north east] at (p0) {$0$}; \node[lab, anchor=north west] at (p2) {$2$};
  \node[lab, anchor=south west] at (p8) {$8$}; \node[lab, anchor=south east] at (p6) {$6$};
  \node[lab, anchor=west, xshift=1.2pt] at (p4) {$4$};
  \node[font=\footnotesize, anchor=north] at (0.8,-0.3) {#5};
 \end{scope}}
\tpanel{0}{0}{1}{1}{$\iota$}
\tpanel{3.8}{0}{0}{1}{$\iota\oplus j_5$}
\tpanel{7.6}{0}{1}{0}{$\iota\oplus j_5\oplus j_2$}
\tpanel{11.4}{0}{0}{0}{$\iota\oplus j_2$}
\draw[cay] (2.05,0.8) -- node[above, font=\footnotesize, black] {$\otimes V_8^{\chi_5}$} (3.35,0.8);
\draw[cay] (5.85,0.8) -- node[above, font=\footnotesize, black] {$\otimes V_8^{\kappa}$} (7.15,0.8);
\draw[cay] (9.65,0.8) -- node[above, font=\footnotesize, black] {$\otimes V_8^{\chi_5}$} (10.95,0.8);
\end{tikzpicture}
\caption{The four induced actions on $\mathcal{O}(PSU(2)_8)$ in Example~\ref{ex:torsor}. The vertex labeled $a$ is $V_a$, filled if it lies in $\mathcal{E}$. Solid arrows show $\widehat{\sigma_3}$, and dashed arrows show $\widehat{\sigma_{11}}$, which is the identity when they are absent. Tensoring with $V_8$ is the half-turn about $V_4$, and the twist $\otimes V_8^\chi$ applies it to $\widehat{\sigma_\ell}$ exactly when $\chi(\sigma_\ell)=-1$. Adjacent panels differ by the indicated twist, and the first and last by $\otimes V_8^\kappa$. Here $\chi_5(\sigma_3)=\kappa(\sigma_3)=-1$, so both twists reverse the $4$-cycle $\widehat{\sigma_3}$. Also $\chi_5(\sigma_{11})=1$ and $\kappa(\sigma_{11})=-1$, so only the twist by $V_8^\kappa$ changes $\widehat{\sigma_{11}}$.}
\label{fig:torsor}
\end{figure}

\begin{example}\label{ex:torsor}
Let $\mathcal{B}:=PSU(2)_8\subset SU(2)_8$, with simple objects $V_0,V_2,V_4,V_6,V_8$ and symmetric center $\mathcal{E}=\{V_0,V_8\}\simeq\mathrm{Rep}(\mathbb{Z}/2)$. Let $\iota$ be the inclusion $\mathcal{B}\subset SU(2)_8$, and let $j_5$ and $j_2$ be the Galois-modular extensions of $\mathcal{E}$ given by $SO(5)_2$ and $\mathrm{Ising}^{\boxtimes2}$ (Example~\ref{ex:generators}). Then $u_{j_5}=V_8^{\chi_5}$ and $u_{j_2}=V_8^{\kappa}$, where $\chi_5(\sigma_\ell)=\bigl(\frac{\ell}{5}\bigr)$ is the character of $\mathbb{Q}(\sqrt5)$. By Theorem~\ref{thm:main}, the four extensions in Figure~\ref{fig:torsor} induce four distinct actions, so no two of them are integrally equivalent. The actions are the orbit of $\hat\sigma_\iota$ under the elements $V_8^\chi$ of $\mathcal{U}(\mathcal{E})$ with $\chi\in\{1,\chi_5,\kappa,\chi_5\kappa\}$, the characters of $\mathbb{Q}$ and of the real quadratic fields of conductor dividing $40$. They factor through $\operatorname{Gal}(\mathbb{Q}(\zeta_{40})/\mathbb{Q})$, which is generated by $\sigma_3$, $\sigma_{11}$ and complex conjugation, the last acting trivially.
\end{example}

\subsection{Canonical actions}\label{subsec:canonical}

In two situations there is a distinguished action. Set
\[
\mathcal{U}_{\mathrm{min}}(\mathcal{E}):=\{u_\mathcal{K}:\mathcal{K}\in\mathrm{Mext}(\mathcal{E})\}\subseteq\mathcal{U}(\mathcal{E}).
\]
Here $u_\mathcal{K}:=u_{j_\mathcal{K}}$ for $j_\mathcal{K}\colon\mathcal{E}\subseteq\mathcal{K}$. It is a subgroup, since $\mathcal{K}\mapsto u_\mathcal{K}$ is a homomorphism on $\mathrm{Mext}(\mathcal{E})$ by part (a) of Theorem~\ref{thm:main}; by Lemma~\ref{lem:symmetric}(c), its kernel consists of the integral extensions (cf.\ \cite{Nikshych2022}).

\begin{proposition}\label{prop:canonical}
\leavevmode
\begin{enumerate}[label=(\alph*)]
\item Every element of $\mathcal{U}_{\min}(\mathcal{E})$ has the form $g^\kappa$ with $g\in\mathcal{E}_{\mathrm{pt}}[2]$. If $\mathcal{E}$ is Tannakian, $\mathcal{U}_{\min}(\mathcal{E})$ is trivial.
\item Suppose $\mathcal{B}$ has a minimal modular extension. We call the actions induced by its minimal modular extensions the \emph{minimal actions} of $\mathcal{B}$. They form one $\mathcal{U}_{\min}(\mathcal{E})$-orbit, so they all agree on $\operatorname{Gal}(\overline{\mathbb{Q}}/\mathbb{Q}(\sqrt2))$. If $\mathcal{E}$ is Tannakian, there is only one (Corollary~\ref{cor:B}), and we denote it by $\hat\sigma^{\min}_\mathcal{B}$.
\item A Galois-modular extension $\iota\colon\mathcal{B}\to\mathcal{M}$ satisfies $\hat\sigma_\iota(\mathbf{1})=\mathbf{1}$ for all $\sigma$ if and only if $\mathcal{M}$ is integral. Such extensions exist if and only if $\mathcal{B}$ is integral, and then they all induce the same action, the \emph{central action} $\hat\sigma^\mathcal{Z}_\mathcal{B}$. One of them is $\mathcal{B}\subseteq\mathcal{Z}(\mathcal{B})$.
\item Let $\mathcal{B}$ be integral with Tannakian symmetric center, and suppose it has a minimal modular extension. Then $\hat\sigma^{\min}_\mathcal{B}=\hat\sigma^\mathcal{Z}_\mathcal{B}$ if and only if its minimal modular extensions are integral.
\end{enumerate}
\end{proposition}

\begin{proof}
For the proof of (a), let $\mathcal{K}\in\mathrm{Mext}(\mathcal{E})$, and let $\mathrm{Rep}(G)\subseteq\mathcal{E}$ be the maximal Tannakian subcategory, with regular algebra $A$. The category $\mathcal{K}_A$ of all $A$-modules is a $G$-crossed braided fusion category, graded by $G$, with trivial component $\mathcal{K}^0_A$ \cite[Section 4.4]{MR2609644}. The trivial component is a minimal modular extension of $\mathcal{E}_G\in\{\mathrm{Vec},\mathrm{sVec}\}$, so it is $\mathrm{Vec}$ or one of the categories $SO(n)_1$, $n\in\mathbb{Z}/16$ \cite{MR3641612}, \cite[Theorem 4.25]{MR3613518}. Comparing $\dim(\mathcal{K}_A)=\dim(\mathcal{E})^2/|G|$ with $\dim(\mathcal{K}_A^0)=\dim(\mathcal{E})^2/|G|^2$ shows that the grading is faithful, so every component has dimension $\dim(\mathcal{K}^0_A)\in\{1,4\}$. The simple objects of $\mathcal{K}=(\mathcal{K}_A)^G$ have Frobenius--Perron dimensions $|O|\cdot\dim\rho\cdot d_Z$, where $d_Z:=\mathrm{FPdim}(Z)$, $Z$ runs over the simple objects of $\mathcal{K}_A$, $O$ is the $G$-orbit of $Z$ and $\rho$ is an irreducible projective representation of its stabilizer \cite[Corollary 2.13]{MR3059899}.

If $\mathcal{E}$ is Tannakian, then $\mathcal{K}_A^0=\mathrm{Vec}$, so every component is a single invertible object. Hence $\mathcal{K}$ is integral, and $u_\mathcal{K}\equiv\mathbf{1}$ by Lemma~\ref{lem:symmetric}(b),(c).

Otherwise, let $Z$ be simple in some component. Then $Z\otimes Z^*\in\mathcal{K}^0_A$ contains $\mathbf{1}$ once.
\begin{itemize}
\item If $\mathcal{K}^0_A$ is pointed, then $d_Z^2$ is the order of the stabilizer of $Z$ in the group of invertible objects, which lies in $\{1,2,4\}$.
\item If $\mathcal{K}^0_A$ is of Ising type, then $d_Z^2=1+a+b\sqrt2$ with $b\in\mathbb{Z}_{\ge0}$ and $a\in\{0,1\}$, since the fermion of $\mathcal{K}^0_A$ is invertible and so its tensor product with $Z$ is simple. The squared dimensions in a component add up to $4$, and comparing the coefficients of $\sqrt2$ shows that $b=0$ throughout.
\end{itemize}
Thus $d_Z\in\{1,\sqrt2,2\}$, and every simple object of $\mathcal{K}$ has Frobenius--Perron dimension, hence dimension, in $\mathbb{Z}\cup\sqrt2\,\mathbb{Z}$.
By Lemma~\ref{lem:symmetric}(b) and~(c), $u_\mathcal{K}$ is trivial on $\operatorname{Gal}(\overline{\mathbb{Q}}/\mathbb{Q}(\sqrt2))$. Hence $u_\mathcal{K}=g^\kappa$, where $g=u_\mathcal{K}(\sigma)$ for any $\sigma$ with $\kappa(\sigma)=-1$.

For the proof of (b), recall that the minimal modular extensions of $\mathcal{B}$ form a $\mathrm{Mext}(\mathcal{E})$-torsor under Baer sums (Theorem~\ref{thm:LKW}(b)). Now apply part (a) of Theorem~\ref{thm:main} and part (a) above.

For the proof of (c), note that $\hat\sigma_\iota(\mathbf{1})=\mathbf{1}$ if and only if $\sigma\circ\varphi_\mathbf{1}=\varphi_\mathbf{1}$, i.e.\ if and only if $\sigma(d_X)=d_X$ for all $X\in\mathcal{O}(\mathcal{M})$. This holds for all $\sigma$ if and only if $\mathcal{M}$ is integral. If $\mathcal{M}$ is integral, so is $\mathcal{B}$. Conversely, if $\mathcal{B}$ is integral, then $\mathcal{B}\subseteq\mathcal{Z}(\mathcal{B})$ has centralizer $\mathcal{B}^{\mathrm{rev}}$ and $\mathcal{Z}(\mathcal{B})$ is integral. If $\hat\sigma_{\iota'}=\hat\sigma_\iota\otimes u$ and both fix $\mathbf{1}$, then $u(\sigma)=\hat\sigma_{\iota'}(\mathbf{1})=\mathbf{1}$.

For the proof of (d), (b) and (c) show that $\hat\sigma^{\min}_\mathcal{B}=\hat\sigma^\mathcal{Z}_\mathcal{B}$ if and only if $\hat\sigma^{\min}_\mathcal{B}(\mathbf{1})=\mathbf{1}$ for all $\sigma$, if and only if one, equivalently every, minimal modular extension is integral.
\end{proof}

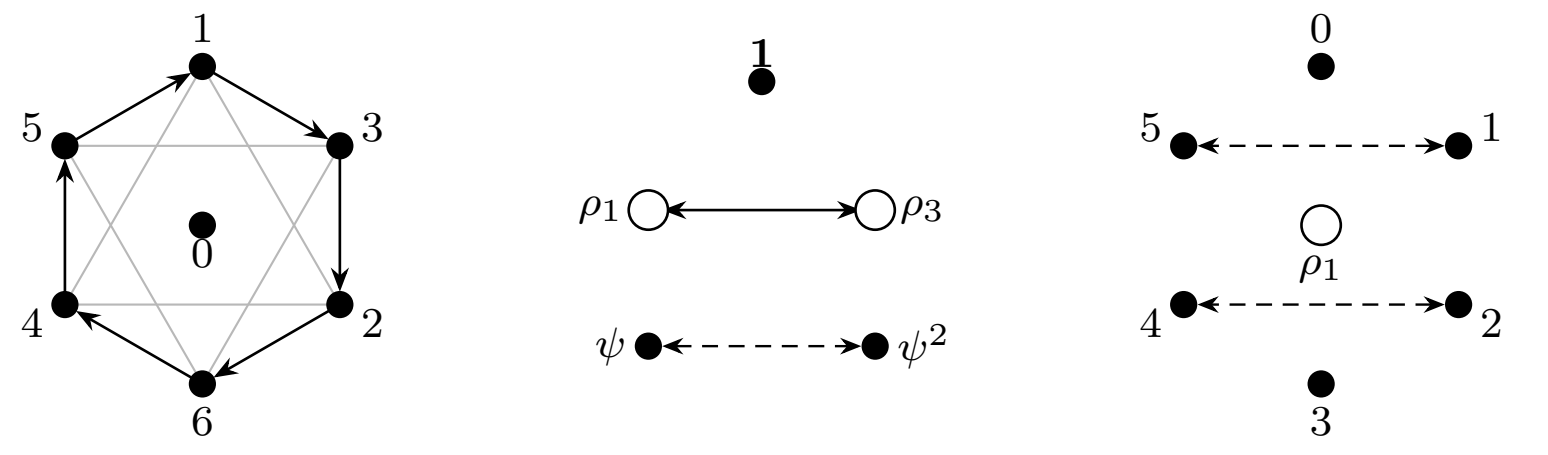
\begin{figure}[ht]
\centering
\begin{tikzpicture}[
  tau/.style={-{Stealth[length=1.6mm,width=1.2mm]}, line width=.5pt, shorten >=2.4pt, shorten <=2.4pt},
  taub/.style={{Stealth[length=1.6mm,width=1.2mm]}-{Stealth[length=1.6mm,width=1.2mm]}, line width=.5pt, shorten >=2.6pt, shorten <=2.6pt},
  eps/.style={{Stealth[length=1.4mm,width=1.1mm]}-{Stealth[length=1.4mm,width=1.1mm]}, densely dashed, line width=.45pt, shorten >=2.6pt, shorten <=2.6pt},
  one/.style={circle, fill=black, inner sep=0pt, minimum size=1.8mm},
  big/.style={circle, fill=white, draw=black, line width=.5pt, inner sep=0pt, minimum size=2.6mm},
  lab/.style={font=\scriptsize, inner sep=1.5pt},
  orbit/.style={black!28, line width=.4pt}]
% (a) Rep(Z/7): chi^t, t in the order of tau (t -> 3t)
\begin{scope}[shift={(0,0)}]
  \foreach \t [count=\i from 0] in {1,3,2,6,4,5} {\coordinate (a\t) at ({90-60*\i}:1.05);}
  \coordinate (a0) at (0,0);
  \draw[orbit] (a1)--(a2)--(a4)--cycle; \draw[orbit] (a3)--(a6)--(a5)--cycle;
  \foreach \s/\t in {1/3,3/2,2/6,6/4,4/5,5/1} {\draw[tau] (a\s)--(a\t);}
  \foreach \t [count=\i from 0] in {1,3,2,6,4,5} {\node[one] at (a\t) {}; \node[lab] at ({90-60*\i}:1.3) {$\t$};}
  \node[one] at (a0) {}; \node[lab, anchor=north, yshift=-1pt] at (a0) {$0$};
  \node[font=\footnotesize] at (0,-1.75) {(a) $\mathrm{Rep}(\mathbb{Z}/7)$};
\end{scope}
% (b) Rep(Z/7 x| Z/3)
\begin{scope}[shift={(3.7,0)}]
  \coordinate (b0) at (0,0.95); \coordinate (r) at (-0.75,0.1); \coordinate (rb) at (0.75,0.1);
  \coordinate (p1) at (-0.75,-0.8); \coordinate (p2) at (0.75,-0.8);
  \draw[taub] (r)--(rb); \draw[eps] (p1)--(p2);
  \node[one] at (b0) {}; \node[lab, anchor=south, yshift=1pt] at (b0) {$\mathbf{1}$};
  \node[big] at (r) {}; \node[lab, anchor=east] at ([xshift=-3pt]r) {$\rho_1$};
  \node[big] at (rb) {}; \node[lab, anchor=west] at ([xshift=3pt]rb) {$\rho_3$};
  \node[one] at (p1) {}; \node[lab, anchor=east] at ([xshift=-2.5pt]p1) {$\psi$};
  \node[one] at (p2) {}; \node[lab, anchor=west] at ([xshift=2.5pt]p2) {$\psi^2$};
  \node[font=\footnotesize] at (0,-1.75) {(b) $\mathrm{Rep}(\mathbb{Z}/7\rtimes\mathbb{Z}/3)$};
\end{scope}
% (c) Rep(Z/7 x| Z/6): psi^k, k = 0..5, and rho_1 of degree 6
\begin{scope}[shift={(7.4,0)}]
  \foreach \k in {0,...,5} {\coordinate (c\k) at ({90-60*\k}:1.05);}
  \draw[eps] (c1)--(c5); \draw[eps] (c2)--(c4);
  \foreach \k in {0,...,5} {\node[one] at (c\k) {}; \node[lab] at ({90-60*\k}:1.3) {$\k$};}
  \node[big] at (0,0) {}; \node[lab, anchor=north, yshift=-4pt] at (0,0) {$\rho_1$};
  \node[font=\footnotesize] at (0,-1.75) {(c) $\mathrm{Rep}(\mathbb{Z}/7\rtimes\mathbb{Z}/6)$};
\end{scope}
\end{tikzpicture}
\caption{The minimal action on $\mathrm{Rep}(\mathbb{Z}/7\rtimes\mathbb{Z}/m)$ in Example~\ref{ex:repG}. Solid arrows show $\hat\tau$ and dashed arrows show $\hat\nu$; each fixes the vertices without an arrow of its kind. Filled vertices have degree $1$. In (a) the vertex $t$ is $\chi^t$, $\hat\tau$ is $t\mapsto3t$, and the grey triangles are the $\mathbb{Z}/3$-orbits $\{1,2,4\}$ and $\{3,5,6\}$. In (b) and (c) the vertex $k$, or $\psi^k$, is the character of $\mathbb{Z}/m$ with $\psi(1)=e^{2\pi i/m}$, inflated to $G_m$.}
\label{fig:repG}
\end{figure}

\begin{example}\label{ex:repG}
Let $\mathcal{B}=\mathrm{Rep}(G)$. Its minimal modular extensions $\mathcal{Z}(\mathrm{Vec}_G^\omega)$ are integral, so $\hat\sigma^{\min}_\mathcal{B}=\hat\sigma^\mathcal{Z}_\mathcal{B}$ by Proposition~\ref{prop:canonical}(d), and this is Galois conjugation $X\mapsto\sigma\circ X$ by the proof of Corollary~\ref{cor:A}. It is the only induced action fixing $\mathbf{1}$, and the only induced action at all if $|G|$ is odd. Figure~\ref{fig:repG} shows it for $G_m:=\mathbb{Z}/7\rtimes\mathbb{Z}/m$, $m\in\{1,3,6\}$, where $\mathbb{Z}/m\subseteq(\mathbb{Z}/7)^\times$ acts by multiplication. All three actions factor through $\operatorname{Gal}(\mathbb{Q}(\zeta_{21})/\mathbb{Q})=\langle\tau\rangle\times\langle\nu\rangle$, where $\tau(\zeta_7)=\zeta_7^3$, $\tau(\zeta_3)=\zeta_3$, $\nu(\zeta_7)=\zeta_7$ and $\nu(\zeta_3)=\zeta_3^{-1}$. The irreducible representations of $G_m$ are the characters $\psi^k$ of $\mathbb{Z}/m$, inflated, and $\rho_t:=\mathrm{Ind}\,\chi^t$ for $t\in(\mathbb{Z}/7)^\times/(\mathbb{Z}/m)$, where $\chi^t(a)=\zeta_7^{ta}$. So the $6$-cycle $\hat\tau$ on $\mathrm{Rep}(\mathbb{Z}/7)$ descends to the orbits: it swaps $\rho_1$ and $\rho_3$ for $m=3$ and fixes $\rho_1$ for $m=6$. The group $G_6$ has the sign character $\mathrm{sgn}:=\psi^3$, so by Corollary~\ref{cor:A} it also has the induced actions $X\mapsto\mathrm{sgn}^\chi(\sigma)\otimes(\sigma\circ X)$, for $\chi$ the character of a real quadratic field; these move $\mathbf{1}$.
\end{example}

\begin{example}\label{ex:canonical}
Let $z:=\psi\boxtimes\psi\in\mathrm{Ising}\boxtimes\mathrm{Ising}^{\mathrm{rev}}$, and let $\mathcal{B}$ be the centralizer of $\{\mathbf{1},z\}$. Its simple objects are $\mathbf{1}$, $z$, $\psi\boxtimes\mathbf{1}$, $\mathbf{1}\boxtimes\psi$ and $\mu\boxtimes\mu$. It is integral, with symmetric center $\{\mathbf{1},z\}\simeq\mathrm{Rep}(\mathbb{Z}/2)$. De-equivariantizing by $z$ gives the toric code, so $\mathcal{B}$ is its equivariantization under the symmetry exchanging $e$ and $m$. Since $\mathrm{Ising}\boxtimes\mathrm{Ising}^{\mathrm{rev}}$ is a minimal modular extension that is not integral, the minimal and central actions differ: $\hat\sigma^{\min}_\mathcal{B}=\hat\sigma^\mathcal{Z}_\mathcal{B}\otimes z^\kappa$. For Galois automorphisms with $\kappa=-1$, such as $\sigma_3$, the minimal action swaps $\mathbf{1}\leftrightarrow z$ and $\psi\boxtimes\mathbf{1}\leftrightarrow\mathbf{1}\boxtimes\psi$, while the central action fixes every simple object.
\end{example}

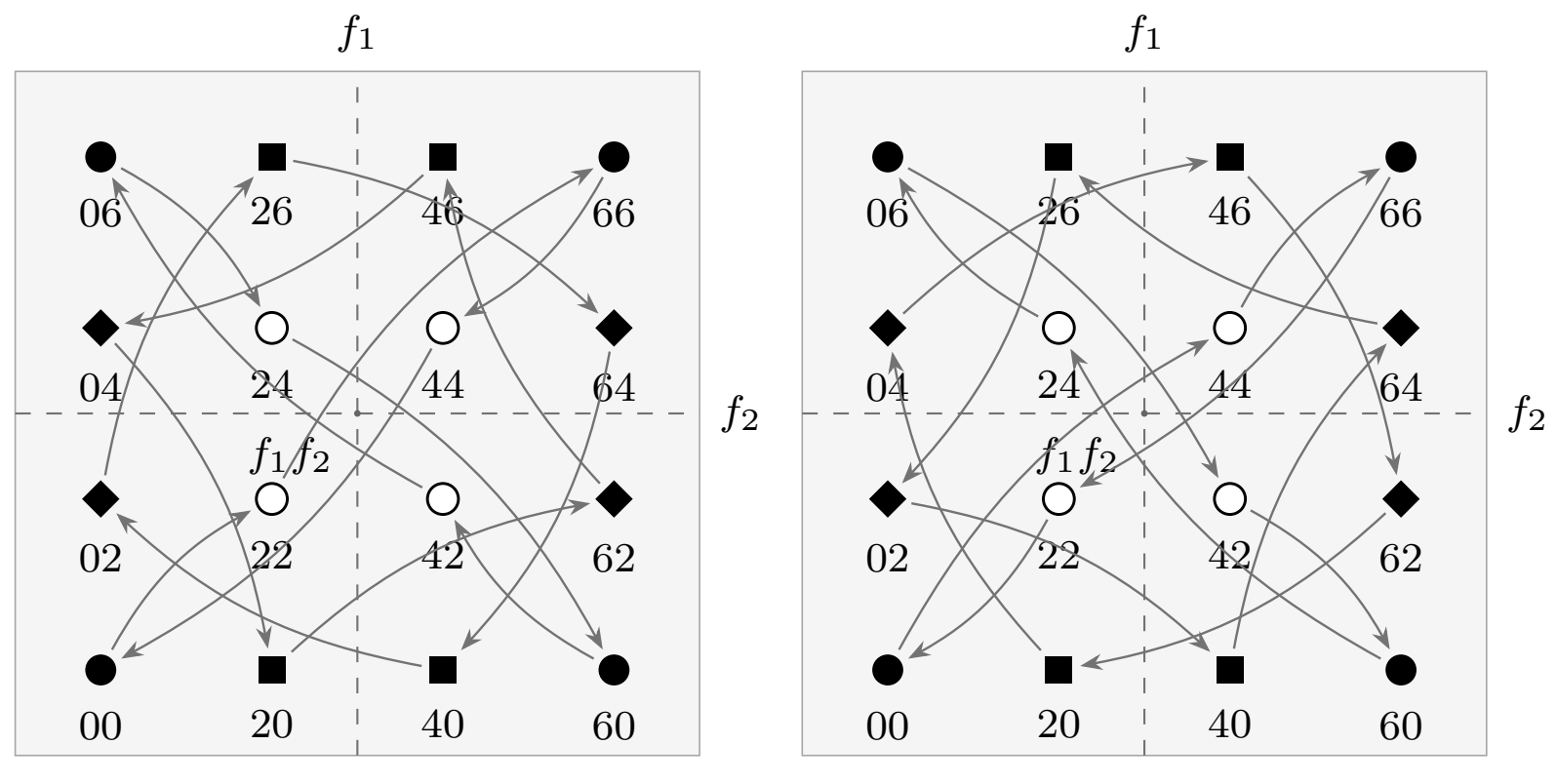
\begin{figure}[ht]
\centering
\begin{tikzpicture}[scale=0.58, gal/.style={-{Stealth[length=1.5mm,width=1.1mm]}, black!55, line width=.45pt, shorten >=1.3pt, shorten <=1.3pt}]
\begin{scope}[xshift=0cm]
\fill[black!4] (0,0) rectangle (8,8);
\draw[black!35, line width=.3pt] (0,0) rectangle (8,8);
\draw[black!55, dashed, line width=.4pt] (4,0)--(4,8) (0,4)--(8,4);
\fill[black!60] (4,4) circle (1.1pt);
\node[font=\scriptsize, anchor=south] at (4,8) {$f_1$};
\node[font=\scriptsize, anchor=west] at (8,4) {$f_2$};
\node[font=\scriptsize, anchor=north east] at (3.95,3.95) {$f_1f_2$};
\node[circle, fill=black, minimum size=2.1mm, inner sep=0pt] (a00) at (1,1) {};
\node[font=\scriptsize, below=0.4mm of a00] {00};
\node[diamond, fill=black, minimum size=2.6mm, inner sep=0pt] (a02) at (1,3) {};
\node[font=\scriptsize, below=0.4mm of a02] {02};
\node[diamond, fill=black, minimum size=2.6mm, inner sep=0pt] (a04) at (1,5) {};
\node[font=\scriptsize, below=0.4mm of a04] {04};
\node[circle, fill=black, minimum size=2.1mm, inner sep=0pt] (a06) at (1,7) {};
\node[font=\scriptsize, below=0.4mm of a06] {06};
\node[rectangle, fill=black, minimum size=1.85mm, inner sep=0pt] (a20) at (3,1) {};
\node[font=\scriptsize, below=0.4mm of a20] {20};
\node[circle, fill=white, draw=black, line width=.55pt, minimum size=2.1mm, inner sep=0pt] (a22) at (3,3) {};
\node[font=\scriptsize, below=0.4mm of a22] {22};
\node[circle, fill=white, draw=black, line width=.55pt, minimum size=2.1mm, inner sep=0pt] (a24) at (3,5) {};
\node[font=\scriptsize, below=0.4mm of a24] {24};
\node[rectangle, fill=black, minimum size=1.85mm, inner sep=0pt] (a26) at (3,7) {};
\node[font=\scriptsize, below=0.4mm of a26] {26};
\node[rectangle, fill=black, minimum size=1.85mm, inner sep=0pt] (a40) at (5,1) {};
\node[font=\scriptsize, below=0.4mm of a40] {40};
\node[circle, fill=white, draw=black, line width=.55pt, minimum size=2.1mm, inner sep=0pt] (a42) at (5,3) {};
\node[font=\scriptsize, below=0.4mm of a42] {42};
\node[circle, fill=white, draw=black, line width=.55pt, minimum size=2.1mm, inner sep=0pt] (a44) at (5,5) {};
\node[font=\scriptsize, below=0.4mm of a44] {44};
\node[rectangle, fill=black, minimum size=1.85mm, inner sep=0pt] (a46) at (5,7) {};
\node[font=\scriptsize, below=0.4mm of a46] {46};
\node[circle, fill=black, minimum size=2.1mm, inner sep=0pt] (a60) at (7,1) {};
\node[font=\scriptsize, below=0.4mm of a60] {60};
\node[diamond, fill=black, minimum size=2.6mm, inner sep=0pt] (a62) at (7,3) {};
\node[font=\scriptsize, below=0.4mm of a62] {62};
\node[diamond, fill=black, minimum size=2.6mm, inner sep=0pt] (a64) at (7,5) {};
\node[font=\scriptsize, below=0.4mm of a64] {64};
\node[circle, fill=black, minimum size=2.1mm, inner sep=0pt] (a66) at (7,7) {};
\node[font=\scriptsize, below=0.4mm of a66] {66};
\draw[gal] (a00) to[bend left=16] (a22);
\draw[gal] (a02) to[bend left=16] (a26);
\draw[gal] (a04) to[bend left=16] (a20);
\draw[gal] (a06) to[bend left=16] (a24);
\draw[gal] (a20) to[bend left=16] (a62);
\draw[gal] (a22) to[bend left=16] (a66);
\draw[gal] (a24) to[bend left=16] (a60);
\draw[gal] (a26) to[bend left=16] (a64);
\draw[gal] (a40) to[bend left=16] (a02);
\draw[gal] (a42) to[bend left=16] (a06);
\draw[gal] (a44) to[bend left=16] (a00);
\draw[gal] (a46) to[bend left=16] (a04);
\draw[gal] (a60) to[bend left=16] (a42);
\draw[gal] (a62) to[bend left=16] (a46);
\draw[gal] (a64) to[bend left=16] (a40);
\draw[gal] (a66) to[bend left=16] (a44);
\node[font=\small, anchor=north] at (4.0,-0.15) {(a) $\hat\sigma_{\mathcal{M}_{00}}$ on $\mathcal{B}$};
\end{scope}
\begin{scope}[xshift=9.2cm]
\fill[black!4] (0,0) rectangle (8,8);
\draw[black!35, line width=.3pt] (0,0) rectangle (8,8);
\draw[black!55, dashed, line width=.4pt] (4,0)--(4,8) (0,4)--(8,4);
\fill[black!60] (4,4) circle (1.1pt);
\node[font=\scriptsize, anchor=south] at (4,8) {$f_1$};
\node[font=\scriptsize, anchor=west] at (8,4) {$f_2$};
\node[font=\scriptsize, anchor=north east] at (3.95,3.95) {$f_1f_2$};
\node[circle, fill=black, minimum size=2.1mm, inner sep=0pt] (b00) at (1,1) {};
\node[font=\scriptsize, below=0.4mm of b00] {00};
\node[diamond, fill=black, minimum size=2.6mm, inner sep=0pt] (b02) at (1,3) {};
\node[font=\scriptsize, below=0.4mm of b02] {02};
\node[diamond, fill=black, minimum size=2.6mm, inner sep=0pt] (b04) at (1,5) {};
\node[font=\scriptsize, below=0.4mm of b04] {04};
\node[circle, fill=black, minimum size=2.1mm, inner sep=0pt] (b06) at (1,7) {};
\node[font=\scriptsize, below=0.4mm of b06] {06};
\node[rectangle, fill=black, minimum size=1.85mm, inner sep=0pt] (b20) at (3,1) {};
\node[font=\scriptsize, below=0.4mm of b20] {20};
\node[circle, fill=white, draw=black, line width=.55pt, minimum size=2.1mm, inner sep=0pt] (b22) at (3,3) {};
\node[font=\scriptsize, below=0.4mm of b22] {22};
\node[circle, fill=white, draw=black, line width=.55pt, minimum size=2.1mm, inner sep=0pt] (b24) at (3,5) {};
\node[font=\scriptsize, below=0.4mm of b24] {24};
\node[rectangle, fill=black, minimum size=1.85mm, inner sep=0pt] (b26) at (3,7) {};
\node[font=\scriptsize, below=0.4mm of b26] {26};
\node[rectangle, fill=black, minimum size=1.85mm, inner sep=0pt] (b40) at (5,1) {};
\node[font=\scriptsize, below=0.4mm of b40] {40};
\node[circle, fill=white, draw=black, line width=.55pt, minimum size=2.1mm, inner sep=0pt] (b42) at (5,3) {};
\node[font=\scriptsize, below=0.4mm of b42] {42};
\node[circle, fill=white, draw=black, line width=.55pt, minimum size=2.1mm, inner sep=0pt] (b44) at (5,5) {};
\node[font=\scriptsize, below=0.4mm of b44] {44};
\node[rectangle, fill=black, minimum size=1.85mm, inner sep=0pt] (b46) at (5,7) {};
\node[font=\scriptsize, below=0.4mm of b46] {46};
\node[circle, fill=black, minimum size=2.1mm, inner sep=0pt] (b60) at (7,1) {};
\node[font=\scriptsize, below=0.4mm of b60] {60};
\node[diamond, fill=black, minimum size=2.6mm, inner sep=0pt] (b62) at (7,3) {};
\node[font=\scriptsize, below=0.4mm of b62] {62};
\node[diamond, fill=black, minimum size=2.6mm, inner sep=0pt] (b64) at (7,5) {};
\node[font=\scriptsize, below=0.4mm of b64] {64};
\node[circle, fill=black, minimum size=2.1mm, inner sep=0pt] (b66) at (7,7) {};
\node[font=\scriptsize, below=0.4mm of b66] {66};
\draw[gal] (b00) to[bend left=16] (b44);
\draw[gal] (b02) to[bend left=16] (b40);
\draw[gal] (b04) to[bend left=16] (b46);
\draw[gal] (b06) to[bend left=16] (b42);
\draw[gal] (b20) to[bend left=16] (b04);
\draw[gal] (b22) to[bend left=16] (b00);
\draw[gal] (b24) to[bend left=16] (b06);
\draw[gal] (b26) to[bend left=16] (b02);
\draw[gal] (b40) to[bend left=16] (b64);
\draw[gal] (b42) to[bend left=16] (b60);
\draw[gal] (b44) to[bend left=16] (b66);
\draw[gal] (b46) to[bend left=16] (b62);
\draw[gal] (b60) to[bend left=16] (b24);
\draw[gal] (b62) to[bend left=16] (b20);
\draw[gal] (b64) to[bend left=16] (b26);
\draw[gal] (b66) to[bend left=16] (b22);
\node[font=\small, anchor=north] at (4.0,-0.15) {(b) $\hat\sigma_{\mathcal{M}_{11}}$ on $\mathcal{B}$};
\end{scope}
\end{tikzpicture}
\caption{The action of $\sigma_3$ in Example~\ref{ex:psu6}, where $\kappa(\sigma_3)=-1$. The object $V_a\boxtimes V_b$ is drawn at $(a+1,b+1)$ and labeled $ab$. The symbol records the dimension: \symdone{} for $1$, \symdh{} and \symdv{} for $1+\sqrt2$ in the first and second factor, and \symdd{} for $(1+\sqrt2)^2$. Dashed lines are the mirrors of $f_1$ and $f_2$, and $f_1f_2$ is the half-turn about the center. Panels (a) and (b) differ by the half-turn $(f_1f_2)^\kappa(\sigma_3)=f_1f_2$.}
\label{fig:psu6}
\end{figure}

\begin{example}\label{ex:psu6}
Let $\mathcal{B}:=PSU(2)_6\boxtimes PSU(2)_6$. Its simple objects are $V_a\boxtimes V_b$ with $a,b\in\{0,2,4,6\}$, and its symmetric center is $\mathcal{E}=\{\mathbf{1},f_1,f_2,f_1f_2\}\simeq\mathrm{sVec}\boxtimes\mathrm{sVec}$. Here $f_1=V_6\boxtimes\mathbf{1}$ and $f_2=\mathbf{1}\boxtimes V_6$ are fermions, and $f_1f_2$ is a boson.

Let $X_0:=SU(2)_6$, and let $X_1$ be the Baer sum of $PSU(2)_6\subset SU(2)_6$ with $\mathrm{sVec}\subset\mathrm{Ising}$, where $c(\mathrm{Ising})=1/2$. So $c(X_0)=9/4$ and $c(X_1)=11/4$, and each $\mathcal{M}_{ij}:=X_i\boxtimes X_j$ is a minimal modular extension of $\mathcal{B}$. By part (a) of Theorem~\ref{thm:main} and Example~\ref{ex:generators}(iii), $\hat\sigma_{\mathcal{M}_{ij}}=\hat\sigma_{\mathcal{M}_{00}}\otimes(f_1^\kappa)^i(f_2^\kappa)^j$. By part (b) of Theorem~\ref{thm:main}, these four actions are distinct, so no two of the $\mathcal{M}_{ij}$ are integrally equivalent. Each $\mathcal{M}_{ij}$ is the Baer sum of $\mathcal{M}_{00}$ with some $\mathcal{K}_{ij}\in\mathrm{Mext}(\mathcal{E})$ (Theorem~\ref{thm:LKW}(b)), and freeness in part (b) of Theorem~\ref{thm:main} gives $u_{\mathcal{K}_{ij}}=(f_1^\kappa)^i(f_2^\kappa)^j$. Together with Proposition~\ref{prop:canonical}(a), this shows $\mathcal{U}_{\min}(\mathcal{E})=\{1,f_1^\kappa,f_2^\kappa,(f_1f_2)^\kappa\}$, and by Proposition~\ref{prop:canonical}(b) the four actions $\hat\sigma_{\mathcal{M}_{ij}}$ are all the minimal actions. The central charges are $9/2$, $11/2$, $5$ and $5$ for $\mathcal{M}_{00}$, $\mathcal{M}_{11}$, $\mathcal{M}_{01}$ and $\mathcal{M}_{10}$, so the central charge does not determine the action.

The picture is a square (Figure~\ref{fig:psu6}). Since $s_{V_a,V_b}$ is proportional to $\sin(\pi(a+1)(b+1)/8)$ in $SU(2)_6$, $\hat\sigma_{\mathcal{M}_{00}}$ multiplies the point $(a+1,b+1)$ by $\ell$ and reflects it back into $[0,8]^2$. The twists $f_1^\kappa$, $f_2^\kappa$ and $(f_1f_2)^\kappa$ act by the two midline reflections and the half-turn.
\end{example}

% ======================================================================
\section{The Galois action on premodular data}\label{sec:data}
% ======================================================================

Throughout this section, $\iota\colon\mathcal{B}\to\mathcal{M}$ is a Galois-modular extension with (additive) central charge $c\in\mathbb{Q}/8\mathbb{Z}$, and $\gamma$ is a cube root of the multiplicative central charge $\xi_\mathcal{M}$. We regard $\mathcal{B}\subseteq\mathcal{M}$, and we write $\hat\sigma$ and $u$ for $\hat\sigma_\iota$ and $u_\iota$. Both depend on $\iota$; the results below hold for every Galois-modular extension $\iota$. By \cite[Theorem 4.1.6]{MR4834521}, $\mathcal{O}(\mathcal{B})$ is $\hat\sigma$-stable, and restricting the defining relation of $\hat\sigma$ to $K(\mathcal{B})$ gives for $Y\in\mathcal{O}(\mathcal{B})$, 
\begin{equation}\label{cond:eq-galois}
\sigma\circ\varphi_Y=\varphi_{\hat\sigma(Y)}\ \text{ on } K(\mathcal{B}),
\end{equation}
where now $\varphi_Y(X)=s_{X,Y}/d_Y$ is computed in $\mathcal{B}$. This is a character of the fusion ring $K(\mathcal{B})$, and $\varphi_{\mathbf{1}}=d$. The \emph{$\mathcal{E}$-orbit} $[Y]$ of $Y\in\mathcal{O}(\mathcal{B})$ is its class under the relation $W\sim Y$ if $W\subset Y\otimes V$ for some $V\in\mathcal{O}(\mathcal{E})$; since $\mathcal{E}$ is closed under tensor products and duals, this is an equivalence relation. Define
\[
K_\mathcal{B}:=\mathbb{Q}(d_X: X\in\mathcal{O}(\mathcal{B})),\quad H_{\mathcal{B}}:=\operatorname{Gal}(\overline{\mathbb{Q}}/K_\mathcal{B}),\quad\text{and}\quad\Phi:=\mathbb{Q}(\varphi_Y(X): X,Y\in\mathcal{O}(\mathcal{B})).
\]
Recall that $\Phi$ lies in a cyclotomic field \cite[Theorem 8.51]{MR2183279}. The \emph{conductor} of a subfield $F$ of a cyclotomic field is the least $f\ge1$ with $F\subseteq\mathbb{Q}(\zeta_f)$. Since $\mathbb{Q}(\zeta_{2m})=\mathbb{Q}(\zeta_m)$ for odd $m$, we have $f\not\equiv2\pmod4$, and $\sigma_\ell$ fixes $F$ whenever $\ell\equiv1\pmod f$. Let $N_\Phi$ and $f_\mathcal{B}$ be the conductors of $\Phi$ and $K_\mathcal{B}$, and define
\[
N_\theta:=\operatorname{lcm}\{\operatorname{ord}(\theta_X): X\in\mathcal{O}(\mathcal{B})\},\qquad\text{and}\qquad e(f):=24f/\gcd(f,12)\quad\text{for}\quad f\geq1.
\]
For modular $\mathcal{B}$, $N_\theta$ is the order of the $T$-matrix and $\Phi$ is the field generated by the ratios $s_{X,Y}/s_{\mathbf{1},Y}$, while the normalized $S$-matrix $s/\sqrt{\dim(\mathcal{B})}$ generates $\Phi(\sqrt{\dim(\mathcal{B})})$. For a root of unity $x$, $x^\ell$ depends only on $\ell$ modulo $\operatorname{ord}(x)$. Since $\widehat{\mathbb{Z}}^\times\to(\mathbb{Z}/m)^\times$ is onto, for every integer $\ell$ prime to $m$ there is some $\sigma_{\hat\ell}$ with $\hat\ell\equiv\ell\pmod m$. We use this to pass between integers $\ell$ and elements of $\operatorname{Gal}(\overline{\mathbb{Q}}/\mathbb{Q})$.

\subsection{Orbits and the unit object}

\begin{lemma}\label{cond:lem-orbits}
Let $\mathcal{B}$ be a pseudounitary braided fusion category and $X,Y,Y',V\in\mathcal{O}(\mathcal{B})$.
\begin{enumerate}[label=(\alph*)]
\item \cite[Proposition 2.5, Corollary 2.14]{MR1990929} $X$ and $V$ centralize each other if and only if $s_{X,V}=d_Xd_V$. In particular $V\in\mathcal{O}(\mathcal{E})$ if and only if $\varphi_V=d$ on $K(\mathcal{B})$.
\item If $V\in\mathcal{O}(\mathcal{E})$ and $W\subset X\otimes V$ is simple, then $\theta_W=\theta_X\theta_V$. Hence $\theta_X^2$ is constant on $\mathcal{E}$-orbits, and $\theta_X$ is constant on $\mathcal{E}$-orbits if $\mathcal{E}$ is Tannakian.
\item $\varphi_Y=\varphi_{Y'}$ on $K(\mathcal{B})$ if and only if $[Y]=[Y']$.
\end{enumerate}
\end{lemma}

\begin{proof}
For the proof of (a), recall that by the balancing axiom $c_{V,X}c_{X,V}$ acts on the $Z$-isotypic part of $X\otimes V$ by the scalar $\omega_Z=\theta_Z/(\theta_X\theta_V)$, a root of unity. Hence $s_{X,V}=\sum_ZN_{XV}^Z d_Z\omega_Z$, while $d_Xd_V=\sum_ZN_{XV}^Zd_Z$. Since every $d_Z>0$, these are equal if and only if every $\omega_Z$ with $N_{XV}^Z\ne0$ equals $1$, i.e.\ if and only if $c_{V,X}c_{X,V}=\mathrm{id}$.

For the proof of (b), note that since $V$ is transparent, $\theta_{X\otimes V}=(\theta_X\otimes\theta_V)c_{V,X}c_{X,V}=\theta_X\otimes\theta_V$, so $\theta_W=\theta_X\theta_V$. The twists of $\mathcal{E}$ are $\pm1$, and they are trivial if $\mathcal{E}$ is Tannakian, so the remaining claims follow.

For the proof of (c), we first show that $\varphi$ is constant on $\mathcal{E}$-orbits. Let $V\in\mathcal{O}(\mathcal{E})$. Since $V$ is transparent, the hexagon axioms identify $c_{Y\otimes V,X}c_{X,Y\otimes V}$ with $(c_{Y,X}c_{X,Y})\otimes\mathrm{id}_V$, and taking traces gives $s_{X,Y\otimes V}=d_Vs_{X,Y}$. Writing $Y\otimes V=\bigoplus_WN_{YV}^WW$, this says
\[
\sum_WN_{YV}^W\,d_W\,\varphi_W=d_Yd_V\,\varphi_Y\qquad\text{on }K(\mathcal{B}).
\]
Distinct characters are linearly independent, and $N_{YV}^Wd_W>0$ for every $W\subset Y\otimes V$. Grouping the left side by character therefore shows that every $W\subset Y\otimes V$ has $\varphi_W=\varphi_Y$.

Next we compute $s^2$. Let $R=\sum_Xd_XX\in K(\mathcal{B})\otimes_\mathbb{Z}\mathbb{R}$. Then $X\cdot R=d_XR$, so a character $\chi$ with $\chi(R)\ne0$ equals $d$. By (a), $\varphi_V(R)=\dim(\mathcal{B})$ for $V\in\mathcal{O}(\mathcal{E})$ and $\varphi_V(R)=0$ otherwise. For fixed $X$, the map $Y\mapsto\varphi_X(Y)=s_{Y,X}/d_X$ is a character of $K(\mathcal{B})$ \cite[Chapter 8]{MR3242743}, so $s_{X,Y}s_{X,W}=d_X\sum_VN_{YW}^Vs_{X,V}$. Hence
\[
(s^2)_{Y,W}=\sum_Xs_{X,Y}s_{X,W}=\sum_VN_{YW}^V\,d_V\,\varphi_V(R)=\dim(\mathcal{B})\sum_{V\in\mathcal{O}(\mathcal{E})}N_{YW}^V\,d_V .
\]
Since $N_{YW}^V=N_{Y^\ast V}^W$, this reads $s^2=\dim(\mathcal{B})\,P\Lambda$. Here $P$ is the permutation matrix of $Y\mapsto Y^\ast$, and $\Lambda_{U,W}=\sum_{V\in\mathcal{O}(\mathcal{E})}d_VN_{UV}^W$.

Finally, we count ranks. The row of $\Lambda$ indexed by $U$ is nonzero, since $\Lambda_{U,U}\ge1$, and it is supported on $[U]$. Rows indexed by representatives of distinct orbits are therefore linearly independent. The column of $s$ indexed by $Y$ is $d_Y\varphi_Y$, so $\operatorname{rank}s$ is at most the number of distinct $\varphi_Y$, hence at most the number of orbits. Thus distinct orbits have distinct characters, as
\[
\#\{\text{orbits}\}\le\operatorname{rank}\Lambda=\operatorname{rank}s^2\le\operatorname{rank}s\le\#\{\varphi_Y\}\le\#\{\text{orbits}\}. \qedhere
\]
\end{proof}

The next lemma describes how $\hat\sigma$ acts on the premodular data of $\mathcal{B}$. Let $D=\sqrt{\dim(\mathcal{B})}>0$ and let $S_{X,Y}=s_{X,Y}/D$ be the normalized $S$-matrix of $\mathcal{B}$.

\begin{lemma}\label{cond:lem-unit}
Let $\sigma=\sigma_\ell$ and $X,Y\in\mathcal{O}(\mathcal{B})$.
\begin{enumerate}[label=(\alph*)]
\item $d_{\hat\sigma(X)}^2=\sigma(d_X)^2\,\dim(\mathcal{B})/\sigma(\dim(\mathcal{B}))$. Hence $\varepsilon_\sigma(X):=\sigma(d_X/D)\,D/d_{\hat\sigma(X)}$ is $\pm1$, and
\[
\sigma(S_{X,Y})=\varepsilon_\sigma(X)\,S_{\hat\sigma(X),Y}=\varepsilon_\sigma(Y)\,S_{X,\hat\sigma(Y)}.
\]
In matrix form, $\sigma(S)=G_\sigma S=SG_\sigma^{T}$, where $G_\sigma$ is the signed permutation matrix with $(G_\sigma)_{X,Y}=\varepsilon_\sigma(X)\,\delta_{\hat\sigma(X),Y}$.
\item $\theta_{\hat\sigma(X)}=\theta_X^{\ell^2}\,\theta_{u(\sigma)}$ and $\theta_{u(\sigma)}=\gamma^{1-\ell^2}$.
\item If $\sigma\in H_{\mathcal{B}}$, then $u(\sigma)\in\mathcal{O}(\mathcal{E})$, so $\theta_{u(\sigma)}=\pm1$, and $\theta_{u(\sigma)}=1$ if $\mathcal{E}$ is Tannakian.
\item If $\sigma\circ\varphi_X=\varphi_X$ on $K(\mathcal{B})$, then $\hat\sigma(X)\in[X]$.
\item $[Y]\mapsto[\hat\sigma(Y)]$ is a well-defined permutation of $\mathcal{E}$-orbits, independent of $\iota$.
\end{enumerate}
\end{lemma}

\begin{proof}
For the proof of (a), note that $\dim(\mathcal{M})=\dim(\mathcal{B})\dim(C_\mathcal{M}(\mathcal{B}))$ with $\dim(C_\mathcal{M}(\mathcal{B}))\in\mathbb{Z}$, so \cite[(4)]{MR4834521} gives
\[
d_{\hat\sigma(X)}^2=\sigma(d_X)^2\,\frac{\dim(\mathcal{M})}{\sigma(\dim(\mathcal{M}))}=\sigma(d_X)^2\,\frac{\dim(\mathcal{B})}{\sigma(\dim(\mathcal{B}))}.
\]
This says exactly that $\varepsilon_\sigma(X)^2=\sigma(d_X^2/\dim(\mathcal{B}))\cdot\dim(\mathcal{B})/d_{\hat\sigma(X)}^2=1$. Now \eqref{cond:eq-galois} gives
\[
\sigma(S_{X,Y})=\sigma\Bigl(\frac{s_{X,Y}}{d_Y}\Bigr)\,\sigma\Bigl(\frac{d_Y}{D}\Bigr)=\frac{s_{X,\hat\sigma(Y)}}{d_{\hat\sigma(Y)}}\cdot\varepsilon_\sigma(Y)\,\frac{d_{\hat\sigma(Y)}}{D}=\varepsilon_\sigma(Y)\,S_{X,\hat\sigma(Y)},
\]
and the other formula follows by symmetry of $S$. The matrix form is a restatement.

For the proof of (b), note that $\sigma^2$ acts on $\theta_X$ and on $\gamma$ by raising to the power $\ell^2$. So the relation $\sigma^2(\gamma^{-1}\theta_X)=\gamma^{-1}\theta_{\hat\sigma(X)}$ of \cite[Theorem II(iii)]{MR3435813} gives $\theta_{\hat\sigma(X)}=\theta_X^{\ell^2}\gamma^{1-\ell^2}$, and the case $X=\mathbf{1}$ gives $\theta_{u(\sigma)}=\gamma^{1-\ell^2}$.

For the proof of (c), note that the values of $\varphi_{\mathbf{1}}=d$ lie in $K_\mathcal{B}$, so $\sigma\circ\varphi_{\mathbf{1}}=\varphi_{\mathbf{1}}$. By \eqref{cond:eq-galois}, $\varphi_{u(\sigma)}=d$ on $K(\mathcal{B})$. So $u(\sigma)\in\mathcal{O}(\mathcal{E})$ by Lemma~\ref{cond:lem-orbits}(a) and we are done since the twists in a Tannakian category are trivial.

For (d), (e), by \eqref{cond:eq-galois} and Lemma~\ref{cond:lem-orbits}(c), $\sigma\circ\varphi_Y$ depends only on $[Y],\sigma$ and determines $[\hat\sigma(Y)]$.
\end{proof}

For modular $\mathcal{B}=\mathcal{M}$, part (a) is the familiar relation $\sigma(S)=G_\sigma S$ \cite{MR1266785}, and part (b) is \cite[Theorem II(iii)]{MR3435813}. For degenerate $\mathcal{B}$ the matrix $S$ is singular, and the relation $\sigma(S)=G_\sigma S$ determines $\hat\sigma(X)$ only up to replacing it by some $W\in[\hat\sigma(X)]$ with $d_W=d_{\hat\sigma(X)}$ (Lemma~\ref{cond:lem-orbits}(c)), e.g.\ $\hat\sigma(X)\otimes g$ with $g\in\mathcal{E}$ invertible. By Theorem~\ref{thm:main}, the permutation $\hat\sigma$ depends on $\iota$ up to a twist by $\mathcal{U}(\mathcal{E})$, but $G_\sigma S=\sigma(S)$ does not depend on $\iota$ at all. For super-modular $\mathcal{B}$, the $\mathcal{E}$-orbits are the pairs $\{X,f\otimes X\}$ labeling the fermionic quotient \cite[Theorem 3.5]{MR3641612}, $\Phi$ is the field $\mathbb{Q}(\hat S)$ of \cite[Remark 4]{BPRZ2021}, the induced permutation of $\mathcal{E}$-orbits is the Galois action of \cite[Section 2.5]{BPRZ2021}, and part (a), written for the unnormalized $S$-matrix, is \cite[Corollary 2.4]{BPRZ2021}.

\subsection{An arithmetic lemma}

\begin{lemma}\label{cond:lem-gcd}
Let $f\ge1$ with $f\not\equiv2\pmod4$, and let $x$ be a root of unity such that $x^{\ell^2-1}=1$ for every integer $\ell\equiv1\pmod f$ prime to $\operatorname{ord}(x)$. Then $\operatorname{ord}(x)\mid e(f)$. Conversely, if $\ell\equiv1\pmod f$ and $p$ is a prime not dividing $\ell$, then $v_p(e(f))\le v_p(\ell^2-1)$. In particular $e(f)\mid\ell^2-1$ for every integer $\ell\equiv1\pmod f$ prime to $6$, so the bound is attained.
\end{lemma}

\begin{proof}
Let $p^a$ exactly divide $\operatorname{ord}(x)$, and write $v_p$ for the $p$-adic valuation. In each case below, choose $\ell$ by the Chinese remainder theorem with the stated property, with $\ell\equiv1\pmod f$, and with $\ell\equiv1\pmod q$ for every other prime $q\mid\operatorname{ord}(x)$. Then $\operatorname{ord}(x)\mid\ell^2-1$, so $a\le v_p(\ell^2-1)$.
\begin{itemize}
\item $p\ge5$, $p\nmid f$: take $\ell\equiv2\pmod p$. Then $\ell^2-1\equiv3\not\equiv0\pmod p$, so $a=0$.
\item $p=3\nmid f$: take $\ell\equiv2\pmod 9$. Then $\ell^2-1\equiv3\pmod9$, so $a\le1$.
\item $p=2\nmid f$: take $\ell\equiv3\pmod{16}$. Then $\ell^2-1\equiv8\pmod{16}$, so $a\le3$.
\item $p$ odd, $v_p(f)=g\ge1$: take $\ell=1+ft$ with $p\nmid t$. Then $v_p(\ell-1)=g$ and $\ell+1\equiv2\pmod p$, so $a\le g$.
\item $p=2\mid f$: then $4\mid f$; let $g=v_2(f)\ge2$ and take $\ell=1+ft$ with $t$ odd. Then $v_2(\ell-1)=g$ and $\ell+1\equiv2\pmod 4$, so $a\le g+1$.
\end{itemize}
On the other hand, $v_p(e(f))=v_p(f)$ for $p\ge5$; $v_3(e(f))=\max(1,v_3(f))$; and $v_2(e(f))=3$ if $f$ is odd, while $v_2(e(f))=v_2(f)+1$ if $4\mid f$. In each case $p^a\mid e(f)$.

For the converse, let $\ell\equiv1\pmod f$ and $p\nmid\ell$. If $p\mid f$, then $v_p(\ell^2-1)\ge v_p(\ell-1)\ge v_p(f)$, and if $p=2$ the factor $\ell+1$ contributes one more $2$. If $p\nmid f$, then $v_p(e(f))$ is $0$ for $p\ge5$, $1$ for $p=3$ and $3$ for $p=2$, and $3\mid\ell^2-1$ for $3\nmid\ell$, $8\mid\ell^2-1$ for odd $\ell$. In each case $v_p(e(f))\le v_p(\ell^2-1)$.
\end{proof}

\subsection{Conductor bounds}\label{sec:conductor}

For a modular fusion category, the entries of the unnormalized $S$-matrix $s$ lie in $\mathbb{Q}(\zeta_N)$ with $N$ the order of $T$ \cite[Proposition 5.7, Section 7]{MR2725181}, and a bound on $N$ in the other direction follows from \cite[Theorem II(iii)]{MR3435813}. The following results extend the second statement to pseudounitary premodular categories, and the first, for $\Phi$, to those with Tannakian symmetric center. For super-modular categories, related constraints were obtained in \cite[Lemma 2.5, Example 3.1]{BPRZ2021} from the $S$- and $T$-matrices of a minimal modular extension; Theorem~\ref{cond:thm}(c) and Corollary~\ref{cond:cor-2group} constrain the twists of $\mathcal{B}$ itself.

\begin{theorem}\label{cond:thm}
Let $\iota\colon\mathcal{B}\to\mathcal{M}$ be a Galois-modular extension with (additive) central charge $c$.
\begin{enumerate}[label=(\alph*)]
\item $N_\theta\cdot c\in\mathbb{Z}$.
\item $\dfrac{f_\mathcal{B}}{\gcd(f_\mathcal{B},12)}\cdot c\in\mathbb{Z}$ if $\mathcal{E}$ is Tannakian, and $\dfrac{2f_\mathcal{B}}{\gcd(f_\mathcal{B},12)}\cdot c\in\mathbb{Z}$ in general.
\item For every $X\in\mathcal{O}(\mathcal{B})$, $\operatorname{ord}(\theta_X)\mid e(N_\Phi)$ if $\mathcal{E}$ is Tannakian, and $\operatorname{ord}(\theta_X^2)\mid e(N_\Phi)$ in general. In particular $N_\theta\mid e(N_\Phi)$ if $\mathcal{E}$ is Tannakian.
\end{enumerate}
\end{theorem}

\begin{proof}
Write $\gamma=\omega\,e^{2\pi ic/24}$ with $\omega^3=1$. In each part we apply Lemma~\ref{cond:lem-gcd} to a root of unity $x$ and a modulus $f$. Given an integer $\ell\equiv1\pmod f$ prime to $\operatorname{ord}(x)$, we use $\sigma=\sigma_{\hat\ell}$ with $\hat\ell\equiv\ell\pmod{\operatorname{lcm}(f,\operatorname{ord}(x))}$. Then $\sigma$ fixes $\mathbb{Q}(\zeta_f)$, and $x^{\hat\ell^2-1}=x^{\ell^2-1}$.

For the proof of (a), take $x=\gamma^{N_\theta}$ and $f=1$. Since $u(\sigma)\in\mathcal{O}(\mathcal{B})$, we have $\theta_{u(\sigma)}^{N_\theta}=1$. Lemma~\ref{cond:lem-unit}(b) then gives $x^{\hat\ell^2-1}=\theta_{u(\sigma)}^{-N_\theta}=1$. By Lemma~\ref{cond:lem-gcd}, $\gamma^{24N_\theta}=1$. Since $\omega^{24}=1$, this says $e^{2\pi i\,cN_\theta}=1$, i.e.\ $N_\theta c\in\mathbb{Z}$.

For the proof of (b), take $f=f_\mathcal{B}$, and $x=\gamma$ if $\mathcal{E}$ is Tannakian, $x=\gamma^2$ in general. Here $\sigma$ fixes $\mathbb{Q}(\zeta_{f_\mathcal{B}})\supseteq K_\mathcal{B}$, so $\sigma\in H_{\mathcal{B}}$. By Lemma~\ref{cond:lem-unit}(b),(c), $\gamma^{1-\hat\ell^2}=\theta_{u(\sigma)}$, which is $1$ if $\mathcal{E}$ is Tannakian and $\pm1$ in general. So $x^{\hat\ell^2-1}=1$. Lemma~\ref{cond:lem-gcd} gives $\gamma^{e(f_\mathcal{B})}=1$, resp.\ $\gamma^{2e(f_\mathcal{B})}=1$. Since $3\mid e(f_\mathcal{B})$, $\omega$ drops out, and $c\,e(f_\mathcal{B})/24=c\,f_\mathcal{B}/\gcd(f_\mathcal{B},12)$.

For the proof of (c), take $f=N_\Phi$, and $x=\theta_X$ if $\mathcal{E}$ is Tannakian, $x=\theta_X^2$ in general. Here $\sigma$ fixes $\Phi\supseteq K_\mathcal{B}$. So $\sigma\in H_{\mathcal{B}}$, and $\sigma\circ\varphi_X=\varphi_X$ on $K(\mathcal{B})$. By Lemma~\ref{cond:lem-unit}(d), $\hat\sigma(X)\in[X]$, and by Lemma~\ref{cond:lem-unit}(b), $\theta_{\hat\sigma(X)}=\theta_X^{\hat\ell^2}\,\theta_{u(\sigma)}$. If $\mathcal{E}$ is Tannakian, then $\theta_{\hat\sigma(X)}=\theta_X$ by Lemma~\ref{cond:lem-orbits}(b) and $\theta_{u(\sigma)}=1$ by Lemma~\ref{cond:lem-unit}(c), so $\theta_X^{\hat\ell^2-1}=1$. In general, $\theta_{\hat\sigma(X)}^2=\theta_X^2$ and $\theta_{u(\sigma)}^2=1$, so squaring gives $(\theta_X^2)^{\hat\ell^2-1}=1$. In both cases $x^{\ell^2-1}=1$, and Lemma~\ref{cond:lem-gcd} applies.
\end{proof}

\vspace{5mm}

\begin{proposition}\label{cond:prop-reverse}
If $\mathcal{E}$ is Tannakian, then $N_\Phi\mid N_\theta$. Hence $N_\Phi\mid N_\theta\mid e(N_\Phi)$.
\end{proposition}

\begin{proof}
Let $\mathcal{E}\simeq\mathrm{Rep}(G)$. The de-equivariantization $\mathcal{B}_G$ is modular, and the free module functor $F\colon\mathcal{B}\to\mathcal{B}_G$ is a braided pivotal functor \cite{MR2609644,MR1936496}, so it preserves traces. Hence $\varphi_Y(X)=s^{\mathcal{B}_G}_{F(X),F(Y)}/d_{F(Y)}$, and this lies in $\mathbb{Q}(\zeta_N)$ with $N=\operatorname{ord}(T_{\mathcal{B}_G})$ \cite[Proposition 5.7, Section 7]{MR2725181}. Every simple object of $\mathcal{B}_G$ is a summand of some $F(X)$ and has twist $\theta_X$, so $N=N_\theta$. The second claim is Theorem~\ref{cond:thm}(c).
\end{proof}

\begin{corollary}\label{cond:cor-2group}
Put $\Phi(\theta)=\Phi(\theta_X:X\in\mathcal{O}(\mathcal{B}))$, $\Phi(\theta^2)=\Phi(\theta_X^2:X\in\mathcal{O}(\mathcal{B}))$ and $L=\operatorname{lcm}(N_\theta,N_\Phi)$.
\begin{enumerate}[label=(\alph*)]
\item If $\mathcal{E}$ is Tannakian, then $\operatorname{Gal}(\Phi(\theta)/\Phi)$ is an elementary abelian $2$-group. In general, the group $\operatorname{Gal}(\Phi(\theta^2)/\Phi)$ is an elementary abelian $2$-group, and $\operatorname{Gal}(\Phi(\theta)/\Phi)$ has exponent dividing $4$.
\item $\operatorname{Gal}(\mathbb{Q}(\zeta_L)/\mathbb{Q}(\zeta_{N_\Phi}))$ is an elementary abelian $2$-group if and only if $N_\theta\mid e(N_\Phi)$. In particular, it is one if $\mathcal{E}$ is Tannakian.
\item If $\mathcal{E}$ is Tannakian, then $\Phi\subseteq\mathbb{Q}(\zeta_{N_\Phi})\subseteq\mathbb{Q}(\zeta_{N_\theta})=\Phi(\theta)$, and both $\operatorname{Gal}(\mathbb{Q}(\zeta_{N_\theta})/\Phi)$ and $\operatorname{Gal}(\mathbb{Q}(\zeta_{N_\theta})/\mathbb{Q}(\zeta_{N_\Phi}))$ are elementary abelian $2$-groups.
\end{enumerate}
\end{corollary}

\begin{proof}
All fields involved lie in a cyclotomic field, so they are abelian over $\mathbb{Q}$. In particular every extension between them is Galois with abelian Galois group, and every element of such a group is the restriction of some $\sigma\in\operatorname{Gal}(\overline{\mathbb{Q}}/\mathbb{Q})$.

For the proof of (a), let $\sigma=\sigma_\ell$ fix $\Phi$. As in the proof of Theorem~\ref{cond:thm}(c), $\sigma\in H_{\mathcal{B}}$, $\hat\sigma(X)\in[X]$ and $\theta_{\hat\sigma(X)}=\sigma^2(\theta_X)\,\theta_{u(\sigma)}$. If $\mathcal{E}$ is Tannakian, $\theta_{\hat\sigma(X)}=\theta_X$ and $\theta_{u(\sigma)}=1$, so $\sigma^2(\theta_X)=\theta_X$. In general both hold up to sign, so $\sigma^2(\theta_X)=\pm\theta_X$. Then $\sigma^2$ fixes $\theta_X^2$, and $\sigma^4$ fixes $\theta_X$. Since $\sigma$ also fixes $\Phi$, the claims follow.

For the proof of (b), note that the group $\operatorname{Gal}(\mathbb{Q}(\zeta_L)/\mathbb{Q}(\zeta_{N_\Phi}))$ consists of the $\sigma_\ell$ with $\ell\in(\mathbb{Z}/L)^\times$ and $\ell\equiv1\pmod{N_\Phi}$. Since $\ell^2\equiv1\pmod{N_\Phi}$ automatically, it has exponent at most $2$ if and only if $N_\theta\mid\ell^2-1$ for every integer $\ell\equiv1\pmod{N_\Phi}$ prime to $N_\theta$. If so, Lemma~\ref{cond:lem-gcd} with $x=\zeta_{N_\theta}$ gives $N_\theta\mid e(N_\Phi)$. Conversely, if $N_\theta\mid e(N_\Phi)$ and $\ell$ is such an integer, then the converse part of Lemma~\ref{cond:lem-gcd} gives $v_p(N_\theta)\le v_p(e(N_\Phi))\le v_p(\ell^2-1)$ for every prime $p\mid N_\theta$. The last claim is Theorem~\ref{cond:thm}(c).

For the proof of (c), note that the $\theta_X$ are roots of unity whose orders have least common multiple $N_\theta$, so they generate $\mathbb{Q}(\zeta_{N_\theta})$. By Proposition~\ref{cond:prop-reverse}, $\mathbb{Q}(\zeta_{N_\Phi})\subseteq\mathbb{Q}(\zeta_{N_\theta})$, so $\Phi(\theta)=\mathbb{Q}(\zeta_{N_\theta})$. Now (a) applies, and $\operatorname{Gal}(\mathbb{Q}(\zeta_{N_\theta})/\mathbb{Q}(\zeta_{N_\Phi}))$ is a subgroup.
\end{proof}

\vspace{2 mm}

\begin{note}\label{cond:note}
\leavevmode
\begin{enumerate}[label=(\roman*),beginpenalty=10000]
\item If $\mathcal{B}$ is modular, then $\mathcal{B}\hookrightarrow\mathcal{B}$ is itself a Galois-modular extension, $N_\theta=\operatorname{ord}(T)$ and $\Phi=\mathbb{Q}(s_{X,Y}/d_Y)$. Theorem~\ref{cond:thm}(a) says that the denominator of $c(\mathcal{B})$ divides $\operatorname{ord}(T)$. For modular $\mathcal{B}$ this is not new: $\xi_\mathcal{B}^2=p^+/p^-$ is a root of unity in $\mathbb{Q}(\zeta_{\operatorname{ord}(T)})$ \cite[Section 7]{MR2725181}, which gives the sharper $\operatorname{lcm}(2,\operatorname{ord}(T))\,c(\mathcal{B})\in4\mathbb{Z}$. The sharper bound fails for degenerate $\mathcal{B}$: for $\mathrm{sVec}\subset\mathrm{Ising}$, $N_\theta=2$ and $c=1/2$, so $\operatorname{lcm}(2,N_\theta)c=1\notin4\mathbb{Z}$. Proposition~\ref{cond:prop-reverse} is the inclusion $\Phi\subseteq\mathbb{Q}(\zeta_{\operatorname{ord}(T)})$ of \cite{MR2725181}, on which its proof relies, and Corollary~\ref{cond:cor-2group}(c) says that this is an elementary abelian $2$-extension, which is implicit in \cite{MR3435813}.
\item Without the Tannakian hypothesis, Proposition~\ref{cond:prop-reverse} fails. For example, $PSU(2)_k$ with $k=6,10,14,18$ has $N_\theta=(k+2)/2$ and $N_\Phi=k+2$ (by direct computation). These examples are consistent with \cite[Conjecture 3.5]{MR3548123}, which predicts $K_\mathcal{B}\subseteq\mathbb{Q}(\zeta_{2N_\theta})$, and show, as Bruillard observes for $PSU(2)_6$, that $2N_\theta$ cannot be replaced by $N_\theta$: there $K_\mathcal{B}=\mathbb{Q}(\sqrt2)\not\subseteq\mathbb{Q}(\zeta_4)$. The group $\operatorname{Gal}(\Phi(\theta)/\Phi)$ need not be elementary abelian, and the squares in Theorem~\ref{cond:thm}(c) are necessary, even for super-modular $\mathcal{B}$. Let $\mathcal{B}$ be the centralizer of $(8,\psi)$ in $(\mathbb{Z}/16,\zeta_{32}^{x^2})\boxtimes\mathrm{Ising}$, where $\theta_\mu=\zeta_{16}$. Then $\Phi=\mathbb{Q}(\zeta_8)$, and $X=(1,\mu)$ has $\theta_X=\zeta_{32}^3$, of order $32\nmid e(8)=48$. Moreover $\sigma_9$ fixes $\Phi$ and $\sigma_9^2(\theta_X)=-\theta_X$, so the exponent in Corollary~\ref{cond:cor-2group}(a) is sharp. The group is elementary abelian if $\mathcal{B}=\mathcal{C}\boxtimes\mathrm{sVec}$ with $\mathcal{C}$ modular: then $\Phi=\mathbb{Q}(s^{\mathcal{C}}_{X,Y}/d_Y)$ and $\theta_{X\boxtimes f}=-\theta_X$, so Corollary~\ref{cond:cor-2group}(a) for $\mathcal{C}$ gives the claim.
\end{enumerate}
\end{note}

\section*{Acknowledgements}
The author acknowledges the use of Claude (Anthropic) in the development of this manuscript and in particular illustrating the induced Galois action in Figures \ref{fig:torsor}--\ref{fig:psu6}. This assistance included suggestions for intermediate claims, proof steps, and corrections, drafting and revision of expository text, literature and citation searches, and computational checks of the examples. All mathematical statements, proofs, computational results, and references were independently checked by the author, who takes sole responsibility for the contents of the paper.

\end{document}